\documentclass[a4paper,10pt,final]{article}

\usepackage[latin1]{inputenc}
\usepackage{amsthm,amsfonts,amsmath}
\usepackage{graphicx,psfrag}
\usepackage{hyperref}
\usepackage{a4wide}
\usepackage{enumerate}

\newtheorem{proposition}{Proposition}[section]
\newtheorem{theorem}[proposition]{Theorem}
\newtheorem{lemma}[proposition]{Lemma}
\newtheorem{corollary}[proposition]{Corollary}
\newtheorem{definition}[proposition]{Definition}

\newenvironment{proofof}[1]{\smallskip\noindent{\textbf{Proof~of~#1.}}%
  \hspace{1pt}}{\hspace{-5pt}{\nobreak\quad\nobreak\hfill\nobreak%
    $\square$\vspace{2pt}\par}\smallskip\goodbreak}

\numberwithin{equation}{section}
\allowdisplaybreaks[4]

\renewcommand{\phi}{\varphi}
\renewcommand{\theta}{\vartheta}
\renewcommand{\epsilon}{\varepsilon}
\renewcommand{\L}[1]{\mathbf{L^#1}}

\newcommand{\C}[1]{\mathbf{C^{#1}}}

\newcommand{\Cc}[1]{\mathbf{C_c^{#1}}}

\newcommand{\BV}{\mathbf{BV}}
\newcommand{\BVloc}{\mathbf{BV_{loc}}}
\newcommand{\modulo}[1]{{\left|#1\right|}}
\newcommand{\norma}[1]{{\left\|#1\right\|}}
\newcommand{\reali}{{\mathbb{R}}}
\newcommand{\naturali}{{\mathbb{N}}}
\newcommand{\Lip}{\mathop\mathbf{Lip}}
\newcommand{\tv}{\mathop\mathrm{TV}}

\newcommand{\pint}[1]{\mathaccent23{#1}}

\renewcommand{\d}[1]{\mathinner{\mathrm{d}{#1}}}

\newcommand{\uno}{J}

\begin{document}

\title{Stability and Optimization\\ in Structured Population Models\\
  on Graphs}

\author{Rinaldo M.~Colombo$^1$ \and Mauro Garavello$^2$}

\footnotetext[1]{INDAM Unit, University of Brescia}

\footnotetext[2]{Department of Mathematics and Applications,
  University of Milano Bicocca}

\maketitle

\begin{abstract}
  \noindent We prove existence and uniqueness of solutions, continuous
  dependence from the initial datum and stability with respect to the
  boundary condition in a class of initial--boundary value problems
  for systems of balance laws. The particular choice of the boundary
  condition allows to comprehend models with very different
  structures. In particular, we consider a juvenile-adult model, the
  problem of the optimal mating ratio and a model for the optimal
  management of biological resources. The stability result obtained
  allows to tackle various optimal management/control problems,
  providing sufficient conditions for the existence of optimal
  choices/controls.

  \medskip


  \medskip

\end{abstract}

\section{Introduction}
\label{sec:I}
This paper is devoted to the following initial--boundary value problem
for a system of balance laws in one space dimension:
\begin{equation}
  \label{eq:1}
  \left\{
    \begin{array}{l@{\quad\qquad}r@{\;}c@{\;}l}
      \displaystyle
      \partial_t u_i
      +
      \partial_x \left(g_i (t,x) \, u_i \right)
      =
      d_i (t,x) \, u_i
      & (t,x) & \in & \reali^+ \times \reali^+
      \\
      \displaystyle
      g_i (t,0) \, u_i (t,0+)
      =
      \mathcal{B}_i \left(t, u_1 (t), \ldots, u_n (t)\right)
      & t & \in & \reali^+
      \\
      \displaystyle
      u_i (0,x) = u_i^o (x)
      & x & \in & \reali^+
    \end{array}
  \right.
  \qquad i=1, \ldots, n\,.
\end{equation}
Here, $i = 1, \ldots, n$ and $t \in \reali^+$ is time. The
\emph{``space''} variable $x$ varies in $\reali^+$ and in the
applications of~\eqref{eq:1} will have the meaning of a
\emph{biological age}, or \emph{size}. The scalar functions $g_1,
\ldots, g_n$ are \emph{growth} functions, $d_1, \ldots, d_n$ are the
death rates and $u_1^o, \ldots, u_n^o$ constitute the initial data. A
key role is played by our choice of the \emph{birth function}
$\mathcal{B}_i$, for $i=1, \ldots, n$, which we assume of the form
\begin{equation}
  \label{eq:B}
  \mathcal{B}_i (t, u_1, \ldots, u_n)
  =
  \alpha_i \! \left(t, u_1 (\bar x_1-), \ldots, u_n (\bar x_n-)\right)
  +
  \beta_i \! \left(
    \int_{I_1} u_1 (x)\d{x}, \ldots, \int_{I_n} u_n (x)\d{x}
  \right)
\end{equation}
for suitable functions $\alpha_i$, $\beta_i$, points $\bar x_i >0$ and
measurable $I_i \subseteq \reali^+$, for $i=1, \ldots, n$.

The literature on equations similar to~\eqref{eq:1} is vast. We refer
for instance to the exhaustive monograph~\cite{PerthameBook} or to the
more recent edition of~\cite{BrauerCastillo-Chavez} and to the
references therein. Specific features of~\eqref{eq:1} are that it is a
\emph{system}, boundary conditions may contain both a \emph{local}
term, the $\alpha_i$, and a \emph{nonlocal} term, the $\beta_i$.

\par From the analytical point of view, in the present treatment we
emphasize the role of the total variation, setting the main result in
$\BV$. In particular, this allows to consider a function of the
type~\eqref{eq:B} and to prove that the boundary data are attained in
the sense of traces, also due to the boundary being non
characteristic. In this setting, the stability of solutions with
respect to $\alpha_i$ and $\beta_i$ is also obtained.  Moreover, the
techniques used in the sequel can easily be extended to more general
source terms as well as to situations where also the space
distribution needs to be taken into account.

\par
\begin{figure}[h!]
  \centering
  \includegraphics[width=0.4\textwidth]{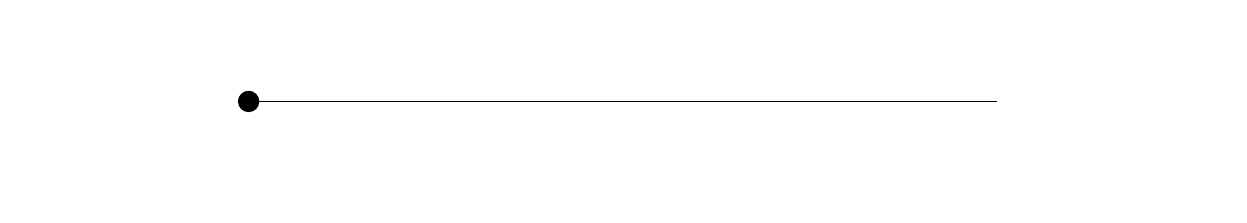}
  \hspace{1.5cm}
  \begin{psfrags}
    \psfrag{j}{$J$} \psfrag{a}{$A$}
    \includegraphics[width=0.4\textwidth]{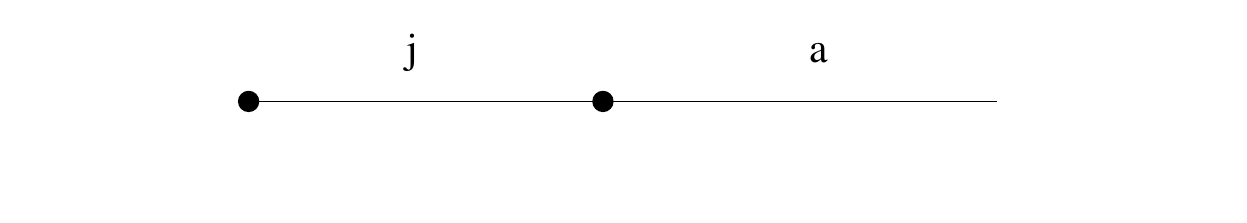}
  \end{psfrags}
  \caption{Biological structures comprised
    in~\eqref{eq:1}--\eqref{eq:B}. Left, a standard \emph{linear}
    setting and, right, a \emph{juvenile-adult} situation.}
  \label{fig:stru1}
\end{figure}

 From the modeling point of view, the use of boundary conditions
of the type~\eqref{eq:B} unifies the treatment of rather diverse
situations. First, it comprises the standard case always covered in
the literature on renewal equations, where the independent variable
$x$ varies along a segment or a half line, see Figure~\ref{fig:stru1},
left.
The dependent variable $u$ represents the population density that at
time $t$ is of size (or age) $x$.

\begin{figure}[h!]
  \centering
  \begin{psfrags}
    \psfrag{m}{$M$} \psfrag{f}{$F$}
    \includegraphics[height=3.4cm]{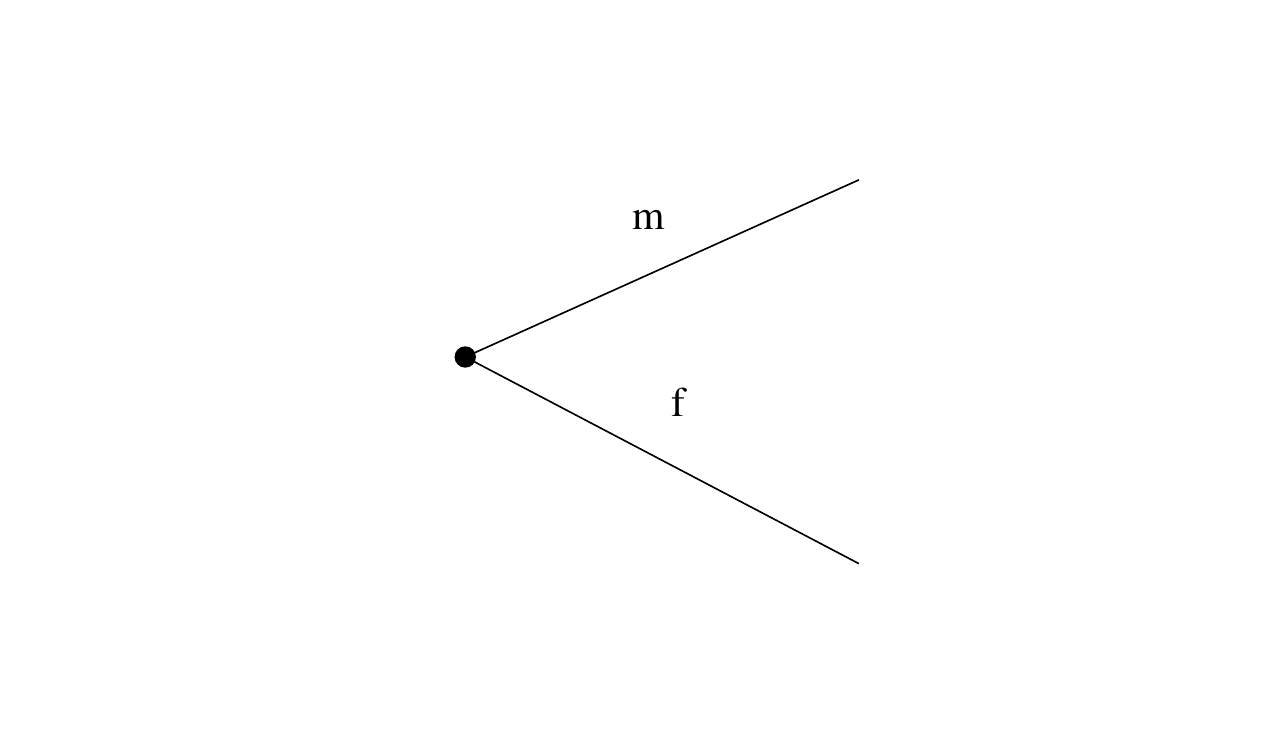}
  \end{psfrags}
  \hspace{1.5cm}
  \begin{psfrags}
    \psfrag{j}{$J$} \psfrag{s}{$S$} \psfrag{r}{$R$}
    \includegraphics[height=3.4cm]{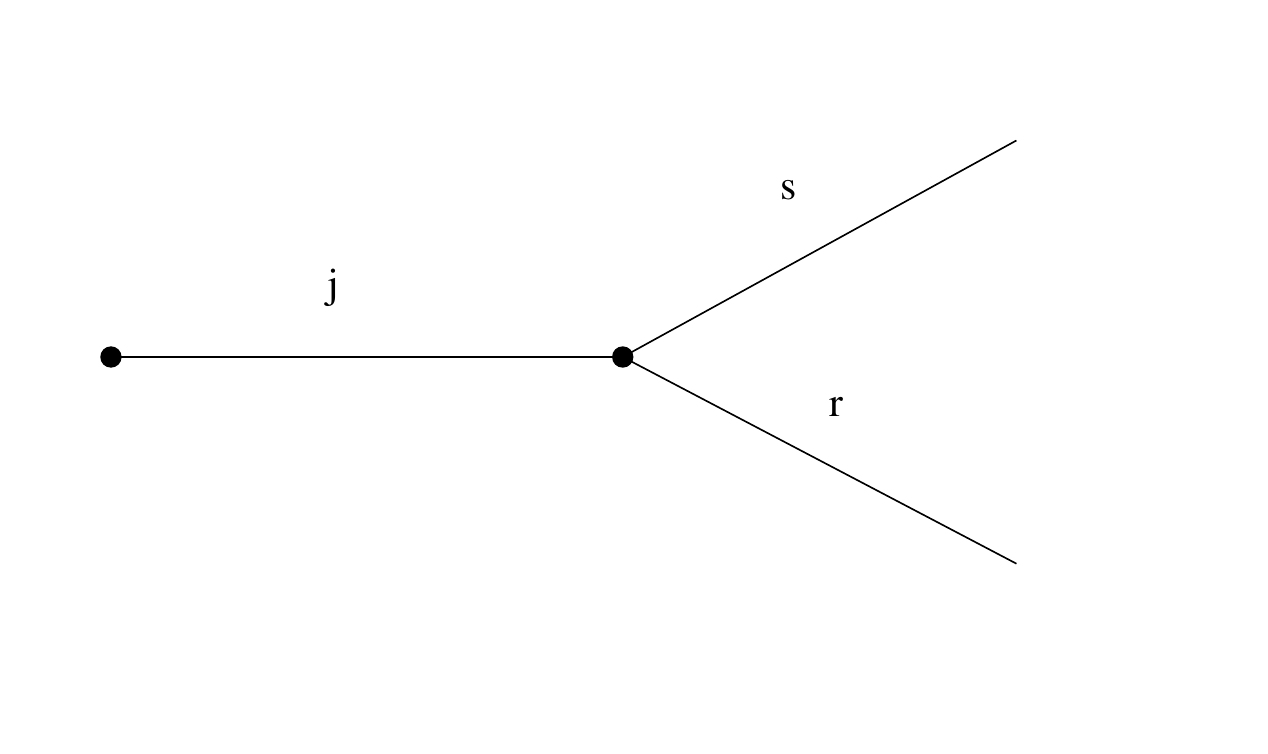}
  \end{psfrags}
  \caption{Biological structures comprised
    in~\eqref{eq:1}--\eqref{eq:B}. Left, a framework corresponding,
    for instance, to sexual reproduction: the two branches correspond
    to males $M$ and females $F$. Right, a structure possibly
    accounting for the exploitation of biological resources: when
    juveniles reach the adult stage, they are split into a part $S$
    which is, say, sold and a part $R$ used for reproduction.}
  \label{fig:stru2}
\end{figure}

A more complicate structure was recently considered
in~\cite{Ackleh2009}, see Figure~\ref{fig:stru1}, right. There, the
size/age biological variable varies along a graph consisting of $2$
distinct sets, corresponding to the juvenile and to the adult stages
in the development of the considered species. Here, we are able to
deal also with this situation, as depicted in Figure~\ref{fig:stru2},
right.

The finite propagation speed intrinsic to models of the
type~\eqref{eq:1} clearly allows to combine various instances of the
graphs above. Other situations of biological interest can be for
instance a three stage linear structure or a tree shaped one, see
Figure~\ref{fig:stru3}. These schemes, as well as many others, all fit
into the scope of Theorem~\ref{thm:wp} below.  In this connection, we
recall that similar network structures are widely considered in the
framework of vehicular traffic modeling,
see~\cite{GaravelloPiccoliBook}.
\begin{figure}[h!]
  \centering
  \includegraphics[height=3.4cm]{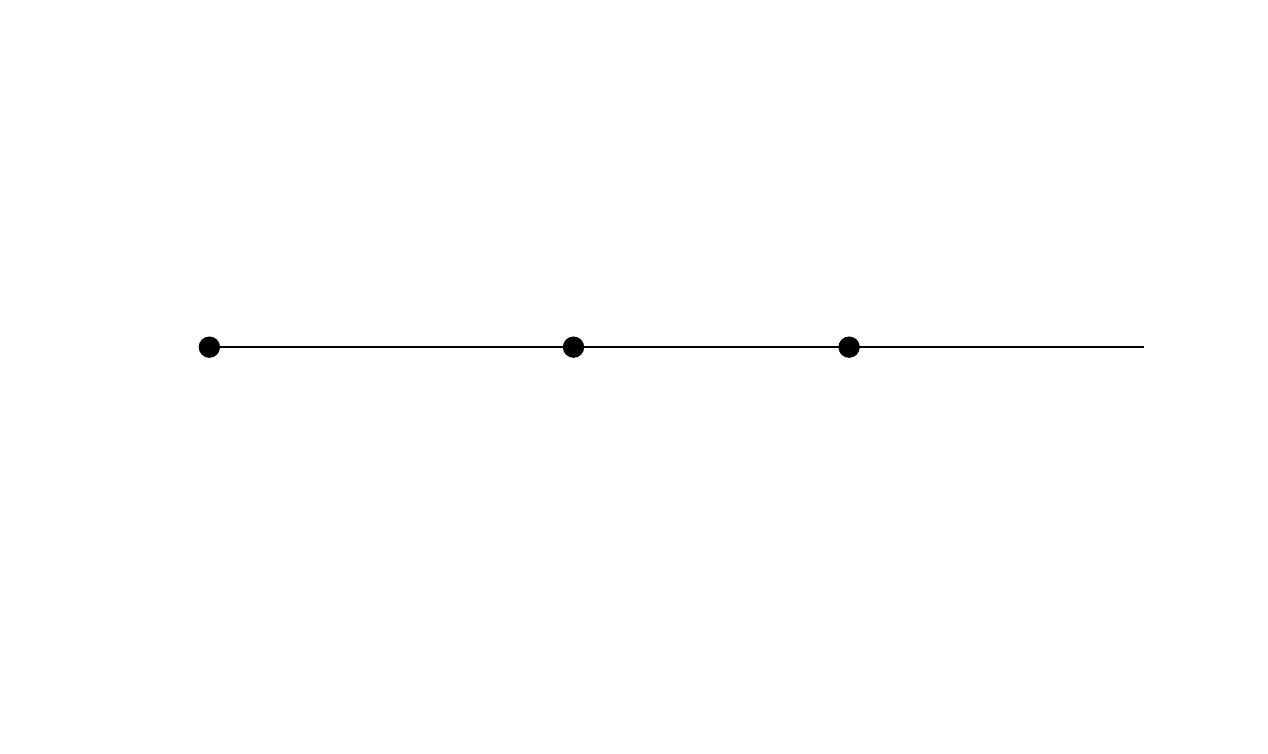}%
  \includegraphics[height=3.4cm]{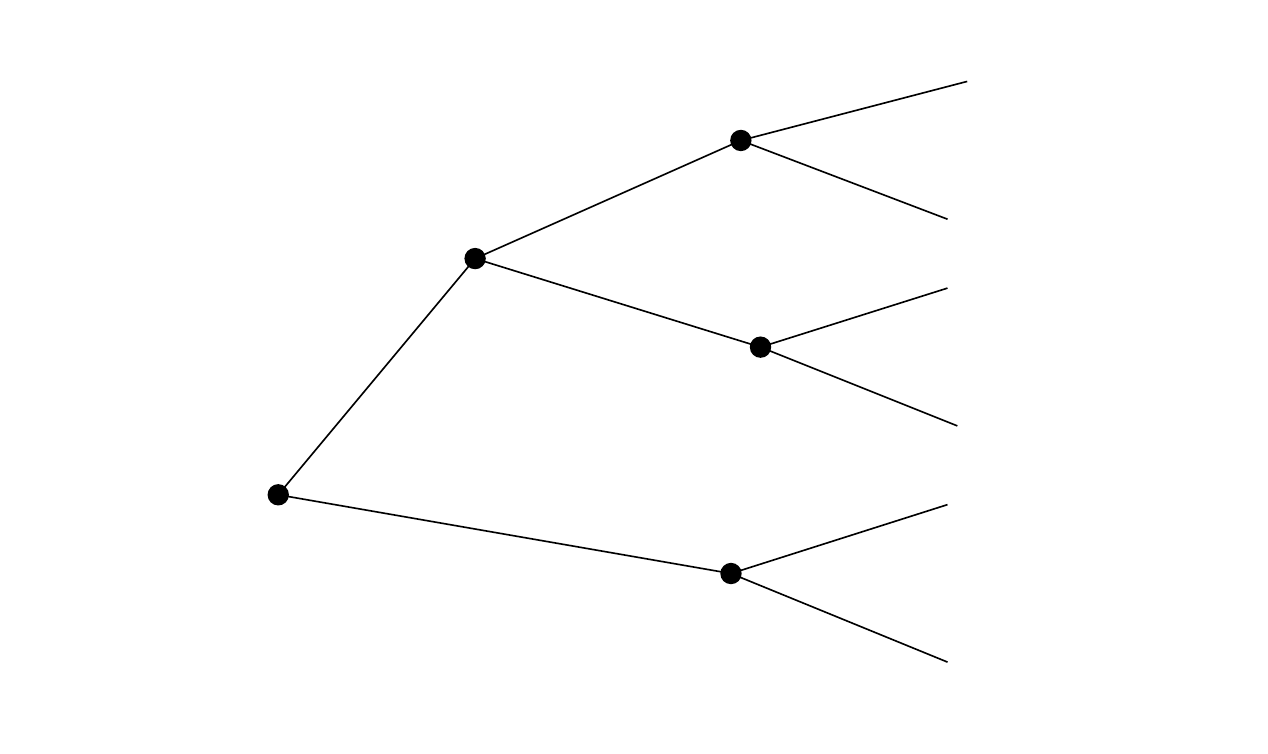}
  \caption{General graphs for further biological structures comprised
    in~\eqref{eq:1}--\eqref{eq:B}.}
  \label{fig:stru3}
\end{figure}

In the case of nonlinear systems of balance laws with flow independent
from the space variable, the initial boundary value problem has been
widely investigated, see for instance~\cite{ColomboGuerra2010}. For
the relations between the problems with boundaries and with junction
see~\cite[Proposition~4.2]{ColomboGuerraHertySchleper}.

The present treatment is self-contained. Section~\ref{sec:AR} is
devoted to the analytical results. Specific applications are in
Section~\ref{sec:Ex}, where sample numerical integrations are also
provided. All technical details are deferred to Section~\ref{sec:TD}.

\section{Analytical Results}
\label{sec:AR}

Throughout, we use the standard notation $\reali^+ = \left[0, +\infty
\right[$ and $\pint{\reali}^+ = \left]0, +\infty\right[$. When $A$ and
$B$ are suitable subsets of $\reali^m$, $\C0 (A;B)$, respectively
$\C{0,1} (A;B)$, $\L1 (A;B)$ or $\L\infty (A;B)$, is the set of
continuous, respectively Lipschitz continuous, Lebesgue integrable or
essentially bounded, maps defined on $A$ and attaining values in
$B$. For the basic theory of $\BV$ functions we refer
to~\cite{AmbrosioFuscoPallara}.

When referring to a function $u \colon \reali^+ \times \reali^+ \to
\reali$, the first argument is time, the second is the biological
age/size variable. If $I \subseteq \reali^+$ is an interval, we denote
\begin{eqnarray*}
  \tv\left(u (t, \cdot)\right)
  & \!\! = \!\! &
  \sup \left\{
    \sum_{h=1}^N \modulo{u (t,x_h) - u (t, x_{h-1})}
    \colon
    N\in \naturali
    \,,\;
    x_{h-1} < x_h, \, \forall h = 1, \ldots, N\!
  \right\}
  \\
  \tv\left(u (t,\cdot);I\right)
  & \!\! = \!\! &
  \sup \left\{
    \sum_{h=1}^N \modulo{u (t,x_h) - u (t, x_{h-1})}
    \colon
    N\in \naturali
    \,,\;
    x_{h-1} < x_h
    \,,\;
    x_h \in I, \, \forall h = 1, \ldots, N \!
  \right\}
  \\
  \tv\left(u (\cdot,x);I\right)
  & \!\! = \!\! &
  \sup \left\{
    \sum_{h=1}^N \modulo{u (t_h,x) - u (t_{h-1}, x)}
    \colon
    N\in \naturali
    \,,\;
    t_{h-1} < t_h
    \,,\;
    t_h \in I, \, \forall h = 1, \ldots, N\!
  \right\} .
\end{eqnarray*}
Preliminarily, we consider the following initial--boundary value
problem for a linear scalar balance law, or \emph{renewal equation}
in~\cite[Chapter~3]{PerthameBook}:
\begin{equation}
  \label{eq:4}
  \left\{
    \begin{array}{l@{\qquad\qquad}r@{\;}c@{\;}l}
      \partial_t u + \partial_x \left(g (t,x) \, u\right) = d (t,x) \,u
      & (t,x) & \in & \reali^+ \times \reali^+
      \\
      u (0, x) = u_o (x)
      & x & \in & \reali^+
      \\
      g (t,0) \, u (t, 0+) = b (t)
      & t & \in & \reali^+
    \end{array}
  \right.
\end{equation}
under the following assumptions
\begin{description}
\item[(b)] $b \in \BVloc \left(\reali^+; \reali\right)$;

\item[(g)] $g \in \C{1} (\reali^+ \times \reali^+; [\check g, \hat
  g])$ for positive $\check g, \hat g$ and $\left\{\begin{array}{r}
      \sup_{t \in \reali^+} \tv\left(g (t,\cdot)\right) < +\infty
      \\
      \sup_{t \in \reali^+} \tv\left(\partial_x g (t,\cdot)\right) <
      +\infty
    \end{array}\right.$;

\item[(d)] $d \in (\C1 \cap \L\infty) (\reali^+ \times \reali^+;
  \reali)$, $\sup_{t \in \reali^+} \tv\left(d (t,\cdot)\right) <
  +\infty$.
\end{description}
\noindent The solutions to~\eqref{eq:4} can be written in terms of the
ordinary differential equation $\dot x = g (t,x)$. If $g$
satisfies~\textbf{(g)}, we can introduce the globally defined maps
\begin{equation}
  \label{eq:XT}
  t \to X (t; t_o, x_o)
  \mbox{ that solves }
  \left\{
    \begin{array}{@{\,}l@{}}
      \dot x = g (t,x)
      \\
      x (t_o) = x_o
    \end{array}
  \right.
  \mbox{ and }
  x \to T (x; t_o, x_o)
  \mbox{ that solves }
  \left\{
    \begin{array}{@{\,}l@{}}
      t' = \frac{1}{g (t,x)}
      \\
      t (x_o) = t_o \,.
    \end{array}
  \right.
\end{equation}
Denote $\gamma (t) = X (t;0,0)$, its inverse being $t = \Gamma
(x)$. Note that
\begin{displaymath}
  \mbox{if } x \geq \gamma (t)
  \mbox{ then } X (0; t,x) \in [0, x]
  \quad \mbox{ and } \quad
  \mbox{if } x < \gamma (t)
  \mbox{ then } T (0; t,x) \in [0,t] \,.
\end{displaymath}
Recall the following definition of solution to~\eqref{eq:4}, see
also~\cite{BardosLerouxNedelec, BressanLectureNotes, Kruzkov,
  PerthameBook, SerreII}.

\begin{definition}
  \label{def:scalar}
  Let~\textbf{(b)}, \textbf{(g)} and~\textbf{(d)} hold. Choose an
  initial datum $u_o \in \L1 (\reali^+; \reali)$. The function $u \in
  \C0 \left(\reali^+; \L1 (\reali^+; \reali)\right)$ is a
  \emph{solution} to~\eqref{eq:4} if
  \begin{enumerate}
  \item for all $\phi \in \Cc1 (\pint{\reali}^+ \times
    \pint{\reali}^+; \reali)$,
    \begin{displaymath}
      \int_{\reali^+} \int_{\reali^+}
      \left(
        u (t,x) \, \partial_t \phi (t,x)
        +
        g (t,x) \, u (t,x) \, \partial_x \phi (t,x)
        +
        d (t,x) \, u (t,x) \, \phi (t,x)
      \right)
      \d{t} \, \d{x} = 0
    \end{displaymath}
  \item $u (0 ) = u_o$
  \item for a.e.~$t \in \reali^+$, $\lim_{x \to 0+} g (t,x) \, u (t,x)
    = b (t)$.
  \end{enumerate}
\end{definition}

\noindent The following Lemma summarizes various properties of the
solution to~\eqref{eq:4}, see also~\cite{PerthameBook}. Here, we
stress the role of $\BV$ estimates. The proof is deferred to
Section~\ref{sec:TD}.

\begin{lemma}
  \label{lem:stability}
  Let~\textbf{(b)}, \textbf{(g)} and~\textbf{(d)} hold. Then, for any
  $u_o \in (\L1 \cap \BV)(\reali^+; \reali)$, the map $u \colon
  \reali^+ \times \reali^+ \to \reali$ defined by
  \begin{equation}
    \label{eq:12}
    \!\!\!
    u (t,x)
    =
    \left\{
      \begin{array}{lr@{\,}c@{\,}l@{}}
        u_o (X (0;t,x)) \,
        \exp\left(
          \int_0^t \left(
            d(\tau,X (\tau;t,x))
            -
            \partial_x g \left(\tau,X (\tau;t,x)\right)
          \right)
          \d\tau
        \right)
        & x & > & \gamma (t)
        \\
        \frac{b\left(T (0;t,x)\right)}{g\left(T (0;t,x),0\right)}
        \exp \left(
          \int_{T (0;t,x)}^t
          \left(
            d(\tau,X (\tau;t,x))
            -
            \partial_x g \left(\tau,X (\tau;t,x)\right)
          \right)
          \d\tau
        \right)
        & x & < & \gamma (t)
      \end{array}
    \right.
  \end{equation}
  solves~(\ref{eq:4}) in the sense of
  Definition~\ref{def:scalar}. Moreover, there exists a constant $C$
  dependent only on $g$ and $d$, see~\eqref{eq:C}, such that the
  following \emph{a priori} estimates hold for all $t \in \reali^+$:
  \begin{eqnarray}
    \label{eq:L1}
    \norma{u (t)}_{\L\infty (\reali^+; \reali)}
    & \leq &
    \left(
      \norma{u_o}_{\L\infty (\reali^+; \reali)}
      +
      \frac{1}{\check g} \, \norma{b}_{\L\infty ([0,t];\reali)}
    \right)
    e^{C t}
    \\
    \label{eq:Linfty}
    \norma{u (t)}_{\L1 (\reali^+; \reali)}
    & \leq &
    \left(
      \norma{u_o}_{\L1 (\reali^+; \reali)}
      +
      \frac{1}{\check g} \, \norma{b}_{\L1 ([0,t];\reali)}
    \right)
    e^{C t}
    \\
    \label{eq:tvx}
    \tv\left(u (t)\right)
    & \leq &
    \left[
      \norma{u_o}_{\L\infty (\reali^+; \reali)}
      +
      \tv (u_o)
      +
      \frac{C+\check g}{\check g^2} \,
      \norma{b}_{\L\infty ([0,t]; \reali)}
      +
      \frac{1}{\check g}
      \tv (b;[0,t])
    \right]
    e^{C t}
    \\
    \label{eq:tvt}
    \!\!\!\!\!\!
    \tv\left(u (\cdot, x); [0,t]\right)
    & \leq &
    \left\{
      \begin{array}{lr@{\,}c@{\,}l}
        \left(
          \norma{u_o}_{\L\infty ([X (0;t,x),x];\reali)}
          +
          \tv (u_o; [X (0;t,x),x])
        \right) e^{C t}
        & x & > &\gamma (t)
        \\[8pt]
        \left(
          \norma{u_o}_{\L\infty ([0,x];\reali)}
          +
          \tv (u_o; [0,x])
        \right) e^{C t}
        \\
        +
        \left(
          \frac{1}{\check g} \tv (b; [0,T (0;t,x)])
          +
          \frac{C\,t}{\check g^2}
          \norma{b}_{\L\infty ([0,T (0;t,x)]; \reali)}
        \right) e^{C t}
        & x & < &\gamma (t) \,.
      \end{array}
    \right.
  \end{eqnarray}
  Moreover, for any interval $I \subseteq \reali^+$,
  \begin{equation}
    \label{eq:tvi}
    \tv \left(\int_I u (\cdot, x) \d{x}; [0,t]\right)
    \leq
    C \int_0^t
    \left(
      \norma{u (\tau)}_{\L\infty (I; \reali)}
      +
      \tv\left(u (\tau,\cdot); I\right)
    \right)
    \d\tau \,.
  \end{equation}
  For every $t \in \reali^+$, there exists a positive $\mathcal L$
  dependent on $\check g, \hat g, C$ and $\tv (b; [0,t])$,
  $\norma{b}_{\L\infty ([0,t]; \reali)}$, such that, for $t', t'' \in
  [0,t]$,
  \begin{equation}
    \label{eq:lip-dependence}
    \norma{u(t') - u(t'')}_{\L1\left(\reali^+; \reali\right)}
    \leq
    \mathcal L \, \modulo{t'' - t'}.
  \end{equation}
  For $u_o', u_o'' \in (\L1 \cap \BV) (\reali^+; \reali)$ and $b',
  b''$ as in~\textbf{(b)}, the solutions $u'$ and $u''$ to
  \begin{equation}
    \label{eq:Two}
    \left\{
      \begin{array}{l}
        \partial_t u + \partial_x \left(g (t,x) \, u\right) = d (t,x) \, u
        \\
        u (0, x) = u_o' (x)
        \\
        g (t,0) \, u (t, 0+) = b' (t)
      \end{array}
    \right.
    \quad \mbox{ and } \quad
    \left\{
      \begin{array}{l@{\qquad\qquad}r@{\;}c@{\;}l}
        \partial_t u + \partial_x \left(g (t,x) \, u\right) = d (t,x) \, u
        \\
        u (0, x) = u_o'' (x)
        \\
        g (t,0) \, u (t, 0+) = b'' (t)
      \end{array}
    \right.
  \end{equation}
  satisfy the stability and monotonicity estimates
  \begin{eqnarray}
    \label{eq:uffa}
    \norma{u' (t)- u'' (t)}_{\L1 (\reali^+; \reali)}
    & \leq &
    \left[
      \norma{u_o'-u_o''}_{\L1 (\reali^+; \reali)}
      +
      \frac{1}{\check g} \, \norma{b'-b''}_{\L1 ([0,t]; \reali)}
    \right]
    e^{C t},
    \\
    \label{eq:mono}
    \left.
      \begin{array}{@{}r@{\;}c@{\;}l@{\mbox{ for all }}r@{\,}c@{\,}l@{}}
        u_o' (x) & \leq & u_o'' (x) & x & \in & \reali^+
        \\
        b' (t) & \leq & b'' (t) & t & \in & \reali^+
      \end{array}
    \right\}
    & \Rightarrow &
    u' (t,x) \leq  u'' (t,x)
    \quad \mbox{ for all } (t,x) \in \reali^+ \times \reali^+ \,.
  \end{eqnarray}
\end{lemma}

Recall that in the present case of a linear conservation law, the
definition of weak solution at~2.~is equivalent to the definition of
Kru\v zkov solution~\cite[Definition~1]{Kruzkov}.

It is immediate to verify that for $u_o = 0$ and $b = 0$,
problem~\eqref{eq:12} admits the solution $u = 0$. Hence, the
monotonicity property~\eqref{eq:mono} also ensures that non-negative
initial and boundary data in~\eqref{eq:1}--\eqref{eq:B} lead to
non-negative solutions.


\smallskip

In order to pass to system~\eqref{eq:1}, we need the following
notation for norms and total variations of functions attaining values
in $\reali^n$:
\begin{displaymath}
  \norma{u}_{\L1 (\reali^+; \reali^n)}
  =
  \sum_{i=1}^n \norma{u_i}_{\L1 (\reali^+; \reali)}
  \,,\quad
  \norma{u}_{\L\infty (\reali^+; \reali^n)}
  =
  \sum_{i=1}^n \norma{u_i}_{\L\infty (\reali^+; \reali)}
  \,,\quad
  \tv (u) = \sum_{i=1}^n \tv (u_i) \,.
\end{displaymath}
As a reference for the usual definition of weak solution to scalar
conservation laws, see~\cite{BardosLerouxNedelec, Kruzkov}.

\begin{definition}
  \label{def:sol}
  Let $T > 0$. Consider~\eqref{eq:1} with $g_1, \ldots, g_n$
  satisfying assumptions~\textbf{(g)} and $d_1, \ldots, d_n$
  satisfying~\textbf{(d)}. Fix an initial datum $u_o \in (\L1 \cap
  \BV) (\reali^+; \reali^n)$. A map
  \begin{displaymath}
    u \in
    \C0\left([0,T]; (\L1 \cap \BV) (\reali^+; \reali^n)\right)
  \end{displaymath}
  is a \emph{solution} to~\eqref{eq:1}--\eqref{eq:B} if, setting
  \begin{displaymath}
    b_i (t)
    =
    \alpha_i\left(t, u_1 (t,\bar x_1-), \ldots, u_n (t, \bar x_n-)\right)
    +
    \beta_i\left(
      \int_{I_1} u_1 (t,x) \d{x}, \ldots, \int_{I_n} u_n (t,x) \d{x}
    \right)
  \end{displaymath}
  for all $i=1 ,\ldots, n$, the $i$-th component $u_i$ is a solution
  to
  \begin{equation}
    \label{eq:questa}
    \left\{
      \begin{array}{l@{\qquad\qquad}r@{\;}c@{\;}l}
        \partial_t u_i + \partial_x \left(g_i (t,x) \, u_i\right)
        =
        d_i (t,x) \, u_i
        & (t,x) & \in & \reali^+ \times \reali^+
        \\
        u_i (0, x) = u_o^i (x)
        & x & \in & \reali^+
        \\
        g_i (t,0) \, u_i (t, 0+) = b_i (t)
        & t & \in & \reali^+
      \end{array}
    \right.
  \end{equation}
  in the sense of Definition~\ref{def:scalar}.
\end{definition}

The following result ensures the well posedness
of~\eqref{eq:1}--\eqref{eq:B}. Its proof is presented in
Section~\ref{sec:TD}.

\begin{theorem}
  \label{thm:wp}
  Let $n \in \naturali \setminus \{0\}$, $\bar x_1, \ldots, \bar x_n
  \in \pint{\reali}^+$, $g_1, \ldots, g_n$ satisfy~\textbf{(g)} and
  $d_1, \ldots, d_n$ satisfy~\textbf{(d)}.  Assume that the maps
  $\alpha \equiv (\alpha_1, \ldots, \alpha_n)$ and $\beta \equiv
  (\beta_1, \ldots, \beta_n)$ satisfy
  \begin{equation}
    \label{eq:5}
    \alpha \in \C{0,1} (\reali^+ \times \reali^n; \reali^n)
    \,,\quad
    \beta \in \C{0,1} (\reali^n; \reali^n)
    \quad \mbox{ and } \quad    \alpha (t,0) = \beta (0) = 0 \,.
  \end{equation}
  Then, for any $u_o \in (\L1 \cap \BV) (\reali^+; \reali^n)$, the
  problem~\eqref{eq:1} admits a unique solution in the sense of
  Definition~\ref{def:sol}. Moreover, there exists an increasing
  function $\mathcal{K} \in \C0 (\reali^+; \reali^+)$ dependent only
  on $\Lip(\alpha)$, $\Lip(\beta)$ and on $C$ in~\eqref{eq:C} such
  that for any initial data $u_o', u_o'' \in (\L1 \cap \BV) (\reali^+;
  \reali^n)$, the corresponding solutions satisfy
  \begin{eqnarray}
    \label{eq:thm1}
    \norma{u' (t) - u'' (t)}_{\L1 (\reali^+; \reali^n)}
    & \leq &
    \mathcal{K} (t)
    \left(
      \norma{u_o'-u_o''}_{\L1 (\reali^+; \reali^m )}
      +
      t \, \norma{u_o'-u_o''}_{\L\infty (\reali^+; \reali^m )}
    \right) \,,
    \\
    \label{eq:thmInfty}
    \norma{u' (t) - u'' (t)}_{\L\infty (\reali^+; \reali^n)}
    & \leq &
    \mathcal{K} (t)
    \left(\norma{u_o'-u_o''}_{\L1 (\reali^+; \reali^m )}
      +
      \norma{u_o'-u_o''}_{\L\infty (\reali^+; \reali^m )}
    \right) \,.
  \end{eqnarray}
  Moreover, if $u_o = 0$, then the solution is $u (t,x) = 0$ for all
  $(t,x) \in \reali^+ \times \reali^+$.
\end{theorem}

We now state separately the stability of solutions
to~\eqref{eq:1}--\eqref{eq:B} with respect to the birth function
$\mathcal{B}$. This result plays a key role in the optimization
problems considered below.

\begin{theorem}
  \label{thm:Stab}
  Let both systems
  \begin{equation}
    \label{eq:7}
    \left\{
      \begin{array}{@{}l@{}}
        \displaystyle
        \partial_t u_i
        +
        \partial_x \left(g_i (t,x) \, u_i \right)
        =
        d_i (t,x) \, u_i
        \\
        \displaystyle
        g_i (t,0) \, u_i (t,0+)
        =
        \mathcal{B}_i' \! \left(t, u_1 (t), \ldots, u_n (t)\right)
        \\
        \displaystyle
        u_i (0,x) = u_i^o (x)
      \end{array}
    \right.
    \quad
    \left\{
      \begin{array}{@{}l@{}}
        \displaystyle
        \partial_t u_i
        +
        \partial_x \left(g_i (t,x) \, u_i \right)
        =
        d_i (t,x) \, u_i
        \\
        \displaystyle
        g_i (t,0) \, u_i (t,0+)
        =
        \mathcal{B}_i'' \! \left(t,u_1 (t), \ldots, u_n (t)\right)
        \\
        \displaystyle
        u_i (0,x) = u_i^o (x)
      \end{array}
    \right.
  \end{equation}
  with
  \begin{eqnarray*}
    \mathcal{B}_i' (t, u_1, \ldots, u_n)
    & = &
    \alpha_i' \! \left(t, u_1 (\bar x_1-), \ldots, u_n (\bar x_n-)\right)
    +
    \beta_i' \! \left(
      \int_{I_1} u_1 (x)\d{x}, \ldots, \int_{I_n} u_n (x)\d{x}
    \right) ,
    \\
    \mathcal{B}_i'' (t, u_1, \ldots, u_n)
    & = &
    \alpha_i'' \! \left(t, u_1 (\bar x_1-), \ldots, u_n (\bar x_n-)\right)
    +
    \beta_i'' \! \left(
      \int_{I_1} u_1 (x)\d{x}, \ldots, \int_{I_n} u_n (x)\d{x}
    \right)
  \end{eqnarray*}
  satisfy the assumptions of Theorem~\ref{thm:wp}. Then, the
  corresponding solutions $u'$ and $u''$ are such that
  \begin{equation}
    \label{eq:stab}
    \norma{u' (t) - u'' (t)}_{\L1 (\reali^+; \reali^n)}
    \leq
    \mathcal{H} (t) \left(
      \norma{\alpha' - \alpha''}_{\C0\left(\reali^+\times\reali^n; \reali^n\right)}
      +
      \norma{\beta' - \beta''}_{\C0\left(\reali^n; \reali^n\right)}
    \right)
  \end{equation}
  where $\mathcal{H} \in \C0 (\reali^+; \reali^+)$ is such that
  $\mathcal{H} (0) = 0$.
\end{theorem}

\noindent The proof is deferred to Section~\ref{sec:TD}.

In applications of Theorem~\ref{thm:wp} to systems motivated by, for
instance, structured population biology, further assumptions are
natural and lead to further reasonable properties.

\begin{proposition}
  \label{prop:bio}
  Under the assumptions of Theorem~\ref{thm:wp}, if the boundary
  functions and the initial data are such that
  \begin{displaymath}
    \begin{array}{rclrcl}
      \partial_{u_j} \alpha_i & \geq & 0
      &
      \mbox{ for all } i,j & = & 1, \ldots, n\,,
      \\
      \partial_{w_j} \beta_i & \geq & 0
      &
      \mbox{ for all } i,j & = & 1, \ldots, n\,,
      \\
      (u_o')_i & \geq & (u_o'')_i
      &
      \mbox{ for all } i & = & 1, \ldots, n \,,
    \end{array}
  \end{displaymath}
  then, the corresponding solutions satisfy $u_i' (t,x) \geq u_i''
  (t,x)$ for all $(t,x) \in \reali^+ \times \reali^+$ and $i=1,
  \ldots, n$.  In particular, if $(u_o')_i \geq 0$ for $i=1, \ldots,
  n$, then $u'_i (t,x) \geq 0$ for $i=1, \ldots, n$.
\end{proposition}

\noindent The proof follows immediately from Theorem~\ref{thm:wp} and
from~\eqref{eq:mono}, hence it is omitted.

\section{Applications}
\label{sec:Ex}

This section is devoted to sample applications of Theorem~\ref{thm:wp}
and Theorem~\ref{thm:Stab} to models inspired by structured population
biology. We selected three cases corresponding to three different
graphs, namely those in Figure~\ref{fig:stru1}, right, and in
Figure~\ref{fig:stru2}.

First, the well posedness ensured by Theorem~\ref{thm:wp} provides a
ground for the reliability of each model. Then, the stability result
in Theorem~\ref{thm:Stab} allows to consider further problems. On the
one hand, it ensures the existence of a choice of parameters in the
equations that lead to solution that best approximate a given set of
data. On the other hand, it allows to tackle the problem of optimal
mating ratio in a population with sexual reproduction. Finally, we
consider the problem of the optimal management of a biological
resource. In the former case, the presentation is based
on~\cite{Ackleh2009, Ackleh2012} where a sensitivity analysis for a
model belonging to the class~\eqref{eq:1}--\eqref{eq:B} is proved. In
the latter cases, we provide numerical integrations showing further
qualitative properties of the models considered.

\subsection{A Nonautonomous Juvenile--Adult Model}
\label{subs:Azmy}

In~\eqref{eq:1}--\eqref{eq:B}, setting $n = 2$ and with reference to
the structure in Figure~\ref{fig:stru1}, right,
\begin{equation}
  \label{eq:AzmyTable}
  \begin{array}{@{}r@{\;}c@{\;}l@{\qquad} r@{\;}c@{\;}l@{\qquad} r@{\;}c@{\;}l@{}}
    u_1 (t,x)& = & J (t,x)
    &
    g_1 (t,x)& = & 1
    &
    \alpha_1 (t,u_1,u_2)& = & 0
    \\
    u_2 (t,x)& = & A (t, x+x_{\min})
    &
    g_2 (t,x)& = & g (t,x+x_{\min})
    &
    \alpha_2 (t,u_1,u_2)& = & u_1
    \\
    \bar x_1 & = & a_{\max}
    &
    d_1 (t,x)& = & -\nu (t,x)
    &
    \beta_1 (w_1,w_2) & = & w_2
    \\
    \bar x_2 & = & 0
    &
    d_2 (t,x)& = & -\mu (t,x+x_{\min})
    &
    \beta_2 (w_1,w_2) & = & 0
  \end{array}
\end{equation}
with moreover $I_2 = [0, x_{\max} - x_{\min}]$, we
recover~\cite[Formula~(2.1)]{Ackleh2009} in the case $\beta=1$, which
we state here for completeness:
\begin{equation}
  \label{eq:Azmy}
  \left\{
    \begin{array}{l@{\qquad\qquad}rcl}
      \partial_t J
      +
      \partial_a J
      =
      -\nu (t,a) \, J
      & (t,a) & \in  & \reali^+ \times [0, a_{\max}]
      \\
      \partial_t A
      +
      \partial_x \left(g (t,x) \, A\right)
      =
      -\mu (t,x) \, A
      &  (t,x) & \in  & \reali^+\times [x_{\min}, x_{\max}]
      \\
      \displaystyle
      J (t,0)
      =
      \int_{x_{\min}}^{x_{\max}} A (t,x) \d{x}
      & t & \in & \reali^+
      \\
      g (t,x_{\min}) \, A (t, x_{\min})
      =
      J (t,a_{\max})
      & t & \in & \reali^+
      \\
      J (0,a) = J_o (a)
      & a & \in & [0, a_{\max}]
      \\
      A (0,x) = A_o (x)
      & x & \in & [x_{\min}, x_{\max}] \,.
    \end{array}
  \right.
\end{equation}
Theorem~\ref{thm:wp} then applies and ensures the well posedness
of~\eqref{eq:Azmy} under assumptions slightly different from those
in~\cite{Ackleh2009}.

\begin{corollary}
  \label{cor:Azmy}
  In~\eqref{eq:Azmy}, assume that
  \begin{displaymath}
    \begin{array}{rclcr@{}}
      \nu & \in & (\C1 \cap \L\infty) (\reali^+ \times [0, a_{\max}]; \reali)
      & \mbox{ and } &
      \sup_{t \in \reali^+} \tv \left(\nu (t, \cdot)\right) < +\infty
      \\
      \mu & \in &
      (\C1 \cap \L\infty) (\reali^+ \times [x_{\min}, x_{\max}]; \reali)
      & \mbox{ and } &
      \sup_{t \in \reali^+} \tv \left(\mu (t, \cdot)\right) < +\infty
      \\
      g & \in & \C1 (\reali^+ \times [x_{\min}, x_{\max}]; [\check g, \hat g])
      & \mbox{ and } &
      \left\{
        \begin{array}{@{}r@{}}
          \sup_{t \in \reali^+} \tv \left(g (t, \cdot)\right) < +\infty
          \\
          \sup_{t \in \reali^+} \tv \left(\partial_x g (t, \cdot)\right) < +\infty
        \end{array}
      \right.
      \\
      J_o & \in & \BV ([0, a_{\max}]; \reali^+)
      \\
      A_o & \in &  \BV ([x_{\min}, x_{\max}]; \reali^+) \,.
    \end{array}
  \end{displaymath}
  Then, problem~\eqref{eq:Azmy} admits a unique solution in the sense
  of Definition~\ref{def:sol}, the continuous dependence
  estimates~\eqref{eq:thm1}--\eqref{eq:thmInfty} and the stability
  estimate~\eqref{eq:stab} apply.
\end{corollary}

For completeness, we remark that the model in~\cite{Ackleh2009}
contains the following slightly more general boundary inflow:
\begin{displaymath}
  J (t,0)
  =
  \int_{x_{\min}}^{x_{\max}} \beta (t,x) \,  A (t,x) \d{x} \,.
\end{displaymath}
As soon as $\beta \in \C1 (\reali^+ \times [x_{\min}, x_{\max}];
[\check \beta , +\infty[)$ for a suitable $\check \beta > 0$, the
change of variables
\begin{equation}
  \label{eq:change}
  \mathcal{A} (t,x) = \beta (t,x) \, A (t,x)
\end{equation}
still allows to apply Theorem~\ref{thm:wp}. Indeed, with this
variable, the second equation in~\eqref{eq:Azmy} becomes
\begin{displaymath}
  \partial_t \mathcal{A}
  +
  \partial_x \left(g (t,x) \, \mathcal{A}\right)
  =
  \left(
    \partial_t \beta (t,x) + g (t,x)\, \partial_x \beta (t,x) - \mu (t,x)
  \right)
  \mathcal{A}\,,
\end{displaymath}
which is again of the type~\eqref{eq:4} and hence Theorem~\ref{thm:wp}
can still be applied.

The stability proved above allows to tackle the problem of parameter
identification. Indeed, through a {Weierstra\ss} argument based on
Theorem~\ref{thm:Stab}, one can prove the existence of a set of
parameters in~\eqref{eq:Azmy} that minimizes a continuous functional
representing the distance between the computed solution and a set of
experimental data. For a detailed sensitivity analysis for a
juvenile--adult model we refer to~\cite{Ackleh2012}.

\subsection{Optimal Mating Ratio}
\label{subs:Mating}

Consider a species consisting of males and females, whose densities at
time $t$ and age $a$ are described through the functions $M = M (t,a)$
and $F = F (t,a)$ on a structure as that in Figure~\ref{fig:stru2},
left. A natural model is then
\begin{equation}
  \label{eq:2Bis}
  \left\{
    \begin{array}{l}
      \displaystyle
      \partial_t M + \partial_a M = - \kappa \, \mu \, M
      \\
      \displaystyle
      \partial_t F + \partial_a F = - (1-\kappa) \, \mu \, F
      \\
      \displaystyle
      M (t,0) + F (t,0)
      =
      \nu
      \min \left\{
        \theta \int_{m_1}^{m_2} M (t,a) \, \d{a}, \;
        (1-\theta) \int_{f_1}^{f_2} F (t,a) \, \d{a}
      \right\}
      \\
      \displaystyle
      \eta \, M (t,0) = (1-\eta) \, F (t,0)
      \\
      \displaystyle
      M (0,a) = M_o (a)
      \\
      F (0,a) = F_o (a)\,.
    \end{array}
  \right.
\end{equation}
Here, $\kappa \, \mu$, respectively $(1-\kappa)\, \mu$, is the
mortality rate of males, respectively females, with $\mu>0$ and
$\kappa \in [0,1]$. The positive parameter $\eta \in [0,1]$ defines
the ratio of male to female newborns, in the sense that every $\eta \,
M$ males, $(1-\eta) \, F$ females are born. The constant $\nu$ is the
fertility rate. We describe the mating ratio at age $a$ through the
parameter $\theta$, with $\theta \in [0,1]$ as follows. The fertile
ages are those in the intervals $[m_1, m_2]$ for males and $[f_1,
f_2]$ for females, where $m_1, m_2, f_1, f_2$ are positive
constants. According to~\eqref{eq:2Bis}, all individuals in their
fertile age might contribute to reproduction provided the condition
imposed by the presence of the mating ratio $\theta$ is met.  If
$\theta \int_{m_1}^{m_2} M (t,a)\d{a}$ exceeds $(1-\theta)
\int_{f_1}^{f_2} F (t,a) \d{a}$, then only $\frac{1-\theta}{\theta} \,
\int_{f_1}^{f_2} F (t,a) \d{a}$ males contribute to the overall
population's fertility.

Problem~\eqref{eq:2Bis} fits into~\eqref{eq:1}--\eqref{eq:B} setting
\begin{displaymath}
  \begin{array}{@{}r@{\;}c@{\;}l@{\qquad} r@{\;}c@{\;}l@{\qquad} r@{\;}c@{\;}l@{}}
    u_1 & = & M
    &
    g_1 & = & 1
    &
    d_1 (t,x)& = & -\kappa
    \\
    u_2 & = & F
    &
    g_2 & = & 1
    &
    d_2 (t,x)& = & - (1-\kappa)\, \mu
    \\
    I_1 & = & [m_1, m_2]
    &
    \alpha_1 & = & 0
    &
    \beta_1 (w_1,w_2) & = & (1-\eta) \, \nu \, \min\{\theta\, w_1, (1-\theta)w_2\}
    \\
    I_2& = & [f_1,f_2]
    &
    \alpha_2 & = & 0
    &
    \beta_2 (w_1,w_2) & = & \eta \, \nu \, \min\{\theta\, w_1, (1-\theta)w_2\}.
  \end{array}
\end{displaymath}

\begin{corollary}
  \label{cor:mating}
  Let $\mu,\nu \in \pint{\reali}^+$; $\eta, \theta, \kappa \in [0,1]$;
  $m_1, m_2, f_1, f_2 \in \reali^+$ with $m_1 < m_2$ and $f_1 <
  f_2$. For $M_o, F_o \in (\L1\cap \BV) (\reali^+; \reali)$,
  problem~\eqref{eq:Azmy} has a unique solution in the sense of
  Definition~\ref{def:sol}, the continuous dependence
  estimates~\eqref{eq:thm1}--\eqref{eq:thmInfty} and the stability
  estimate~\eqref{eq:stab} apply.
\end{corollary}

\noindent The proof is immediate and, hence, omitted. Here, we note
that the presence of $\C1$ positive weights in the integrands defining
the boundary data can be recovered through a change of variables
entirely similar to that in~\eqref{eq:change}.

A first immediate property of the solutions to~\eqref{eq:2Bis} is that
a zero initial density in either of the two sexes leads to the
extinction of the other at exponential speed.

Several different optimization problems can be tackled in the
framework of~\eqref{eq:2Bis}. It is possible to investigate the
relations between the parameters $\kappa$ (identifying relative
mortality), $\eta$ (the relative natality) and $\theta$ (the mating
ratio). Below, we look for the optimal mating ratio for given relative
natality and mortality coefficients.

To this aim, consider the instantaneous average fertility rate over
the fertile population
\begin{equation}
  \label{eq:R}
  R
  =
  \frac{\nu \,
    \min \left\{
      \theta \int_{m_1}^{m_2} M (t,a) \, \d{a}, \;
      (1-\theta) \int_{f_1}^{f_2} F (t,a) \, \d{a}
    \right\}}{\int_{m_1}^{m_2} M (t,a) \, \d{a} +
    \int_{f_1}^{f_2} F (t,a) \, \d{a}} \,.
\end{equation}
Remark that the functions $M$ and $F$ in~\eqref{eq:R} are solutions
to~\eqref{eq:2Bis}, hence they depend on the mating ratio $\theta$
that enters the boundary condition throughout the time interval
$[0,t]$. It is natural to assume that a key role is played by the
maximal value of $R$, which is obtained by the choice
\begin{equation}
  \label{eq:theta}
  \theta
  =
  \frac{\int_{f_1}^{f_2} F (t,a) \, \d{a}}{\int_{m_1}^{m_2} M (t,a) \, \d{a} + \int_{f_1}^{f_2}F (t,a)\, \d{a}} \,.
\end{equation}
Remarkably, this leads to the maximal fertility rate
\begin{displaymath}
  R
  =
  \nu \,
  \frac{\left(\int_{m_1}^{m_2} M (t,a) \, \d{a}\right) \; \left(\int_{f_1}^{f_2} F (t,a) \, \d{a}\right)}{\left(\int_{m_1}^{m_2} M (t,a) \, \d{a} + \int_{f_1}^{f_2} F (t,a) \, \d{a}\right)^2}
\end{displaymath}
coherently with the classical \emph{harmonic mean} law,
see~\cite{Keyfitz1972, Manfredi1993, Schoen1981, Schoen1983,
  Schoen1988, SundelofAberg2006}.

On the other hand, the right hand sides in~\eqref{eq:R}
and~\eqref{eq:theta} are time dependent and it can be hardly accepted
that $\theta$ is instantaneously adjusted to the value that maximizes
$R$. More reasonably, one may imagine that $\theta$ is
optimal\footnote{Here and in the sequel, \emph{optimal} is understood
  in the sense that it is the value that maximizes $R$. However, an
  excessive natality rate might turn out to be not \emph{optimal} from
  the biological point of view.} over a suitably long time interval.
We are thus lead to introduce the utility function
\begin{equation}
  \label{eq:utility-mating}
  \mathcal{R} (\theta; T, M_o,F_o)
  =
  \frac{1}{T}
  \int_0^T
  \frac{\nu \,
    \min \left\{
      \theta \int_{m_1}^{m_2} M (t,a) \, \d{a}, \;
      (1-\theta) \int_{f_1}^{f_2} F (t,a) \, \d{a}
    \right\}}{\int_{m_1}^{m_2} M (t,a) \, \d{a} +
    \int_{f_1}^{f_2} F (t,a) \, \d{a}}
  \d{t} \,.
\end{equation}
We thus consider the problem
\begin{displaymath}
  \mbox{find } \theta
  \mbox{ that maximizes }
  \mathcal{R} (\theta; T, M_o,F_o) \,.
\end{displaymath}
A straightforward corollary of Theorem~\ref{thm:Stab} ensures the
existence of one such $\theta$.

\begin{corollary}
  \label{cor:mating2}
  Under the assumptions of Corollary~\ref{cor:mating}, for any $T \in
  \pint{\reali}^+$ and any initial datum $(M_o,F_o) \in (\L1 \cap \BV)
  (\reali^+; \reali)$ there exists a $\theta_* \in \left]0,1\right[$
  such that
  \begin{displaymath}
    \mathcal{R} (\theta_*; T, M_o,F_o)
    =
    \max_{\theta \in [0,1]} \mathcal{R} (\theta; T, M_o,F_o) \,.
  \end{displaymath}
\end{corollary}

The proof is immediate: thanks to Theorem~\ref{thm:Stab}, the function
$\theta \to \mathcal{R} (\theta; T, M_o,F_o)$ is continuous for any
choice of $T \in \reali^+$ and $(M_o,F_o) \in (\L1 \cap \BV)
(\reali^+; \reali)$. By the compactness of $[0,1]$, {Weierstra\ss}
Theorem ensures the existence of $\theta_*$. Moreover, since
$\mathcal{R} (0; T, M_o,F_o) = \mathcal{R} (1; T, M_o,F_o) = 0$, we
also have $\theta_* \in \left]0,1\right[$.

It can be of interest to note that $M$ and $F$ may well increase
exponentially with time, but $\mathcal{R} (\theta; T, M_o,F_o) \in [0,
\nu]$ for all $T \in \pint{\reali}^+$ and $(M_o,F_o) \in (\L1 \cap
\BV) (\reali^+; \reali)$ .

As a specific example, we consider the situation identified by the
following choices of functions and parameters
in~\eqref{eq:2Bis}--(\ref{eq:utility-mating}):
\begin{equation}
  \label{eq:parameters-mating}
  \begin{array}{@{}rcl@{\qquad}rcl@{\qquad}rcl@{\qquad}rcl@{}}
    \kappa & = & 0.600
    &
    \mu & = & 0.020
    &
    m_1 & = & 18
    &
    f_1 & = & 16
    \\
    \eta & = & 0.485
    &
    \nu & = & 3
    &
    m_2 & = & 60
    &
    f_2 & = & 55
  \end{array}
\end{equation}
and we consider $\theta$ as a control parameter in $[0,1]$.  As
initial datum we choose
\begin{equation}
  \label{eq:idMF}
  M_o (a) = 10
  \quad \mbox{ and } \quad
  F_o (a) = 10
  \quad \mbox{ for all } a \,.
\end{equation}
\begin{figure}[h!]
  \centering
  \includegraphics[width=0.5\textwidth]{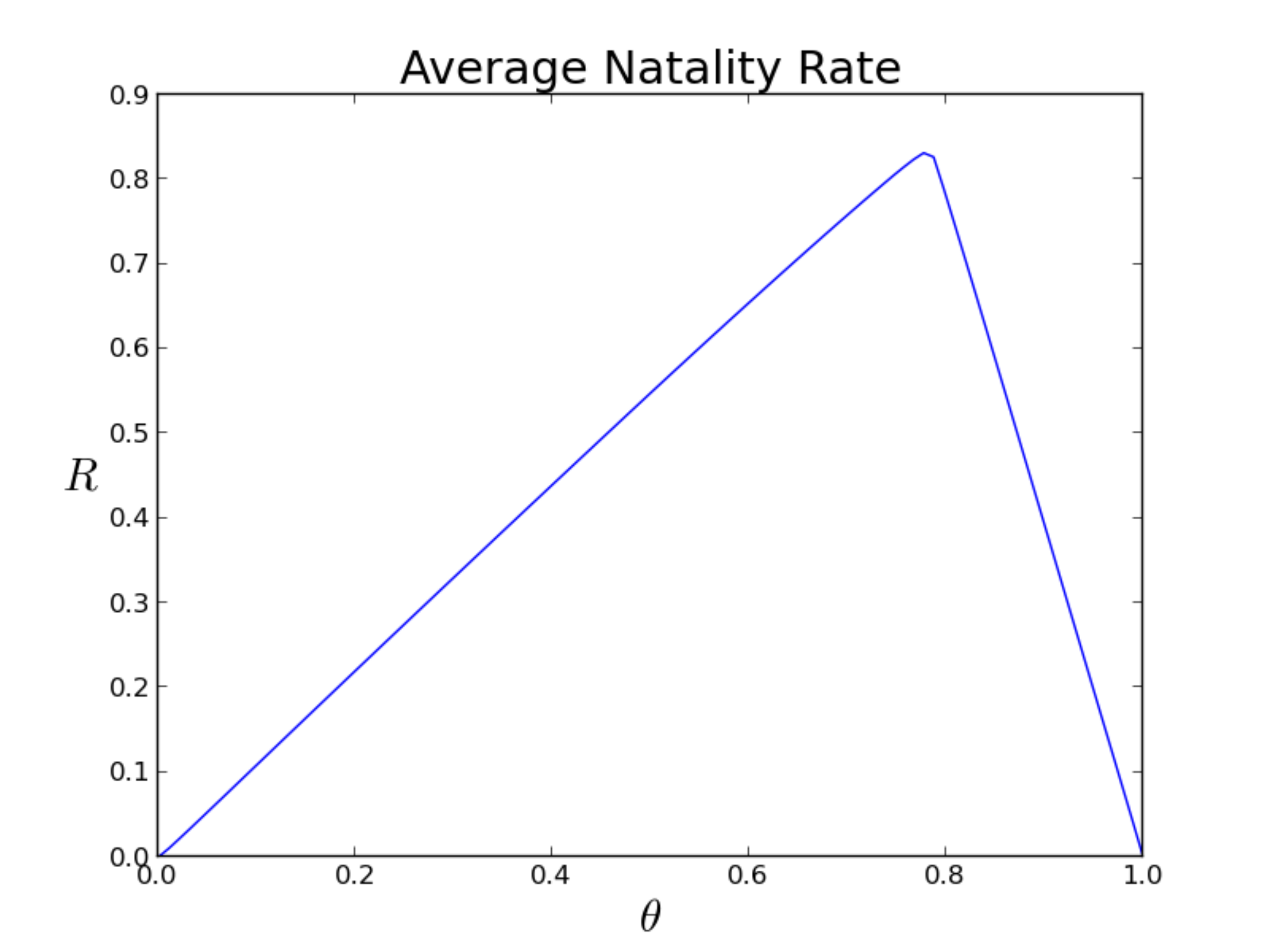}
  \caption{Average fertility rate~(\ref{eq:utility-mating}) along the
    solutions to~\eqref{eq:2Bis}--(\ref{eq:parameters-mating}) with
    initial datum~(\ref{eq:idMF}) plotted as a function of
    $\theta$. For $\theta=0$ or $\theta = 1$, there is no reproduction
    and the population extinguishes.  For $\theta \approx 0.77$, the
    average natality rate reaches its maximum value, which is
    approximately equal to $0.83$.}
  \label{fig:CostMF}
\end{figure}
The graph of the average fertility rate $\mathcal R(\theta; T, M_o,
F_o)$ 
as a function of $\theta$ for $T=500$ is in Figure~\ref{fig:CostMF}.
The outcome shows a reasonable qualitative behavior.  As $\theta \to
0$ or $\theta \to 1$, the number of newborns goes to $0$; hence the
population extinguish.  Near to $\theta = 0.77$ there is an optimal
choice for the parameter $\theta$ with respect to the average
fertility rate~(\ref{eq:utility-mating}), which yields a maximal value
of $0.83$, see~Figure~\ref{fig:contourMF}
\begin{figure}[h!]
  \centering
  \includegraphics[width=0.3\textwidth]{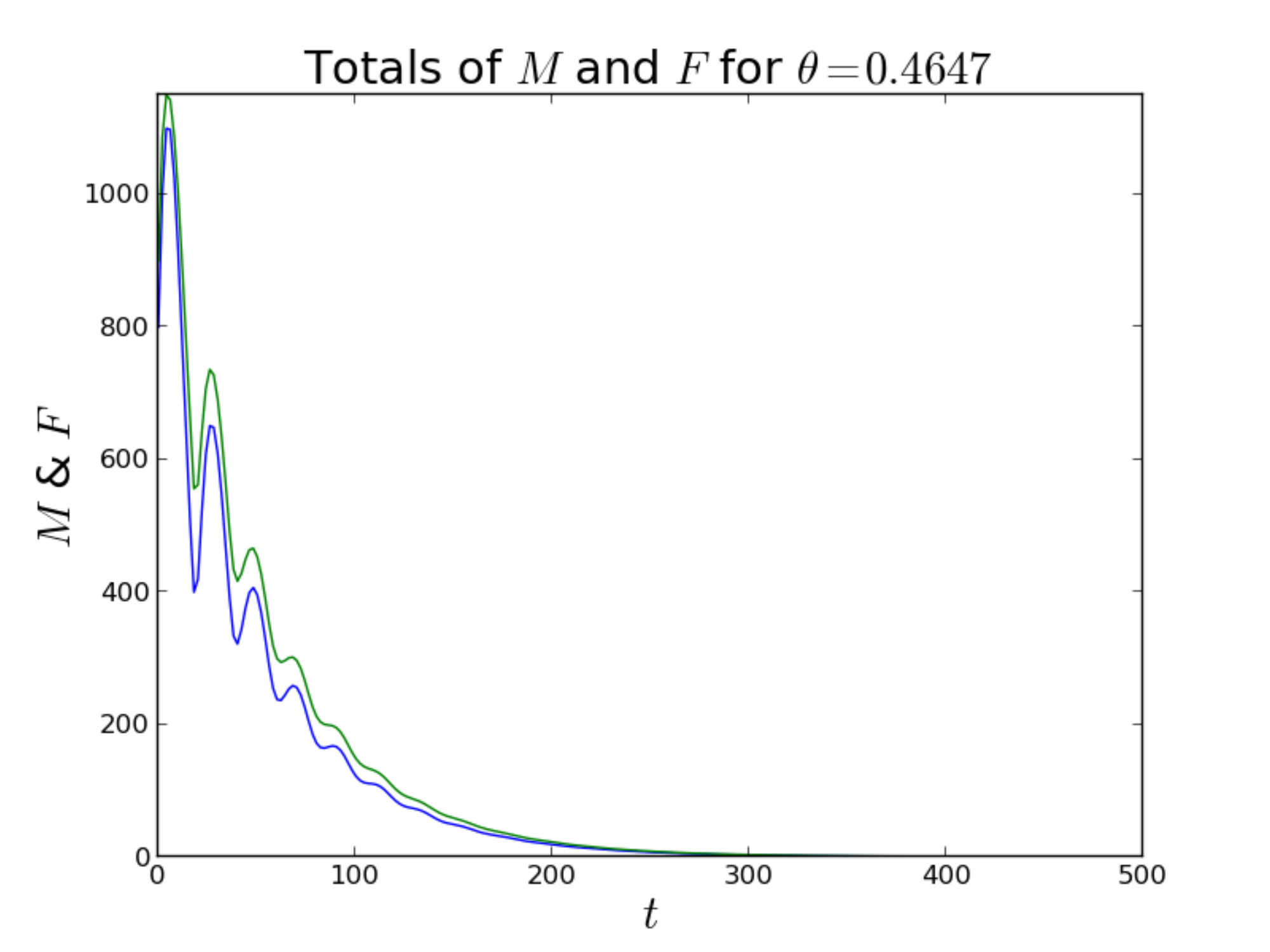}%
  \includegraphics[width=0.3\textwidth]{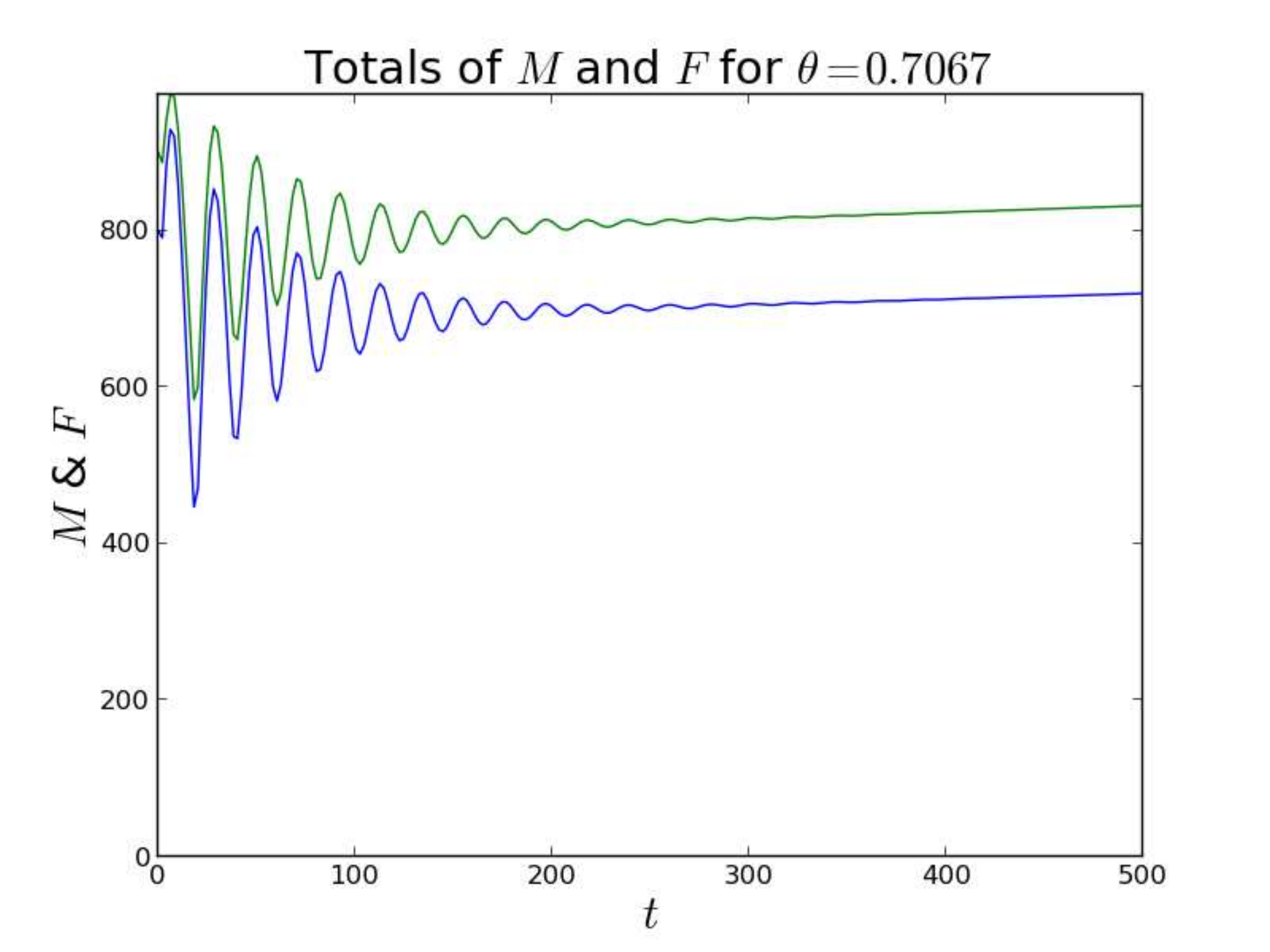}%
  \includegraphics[width=0.3\textwidth]{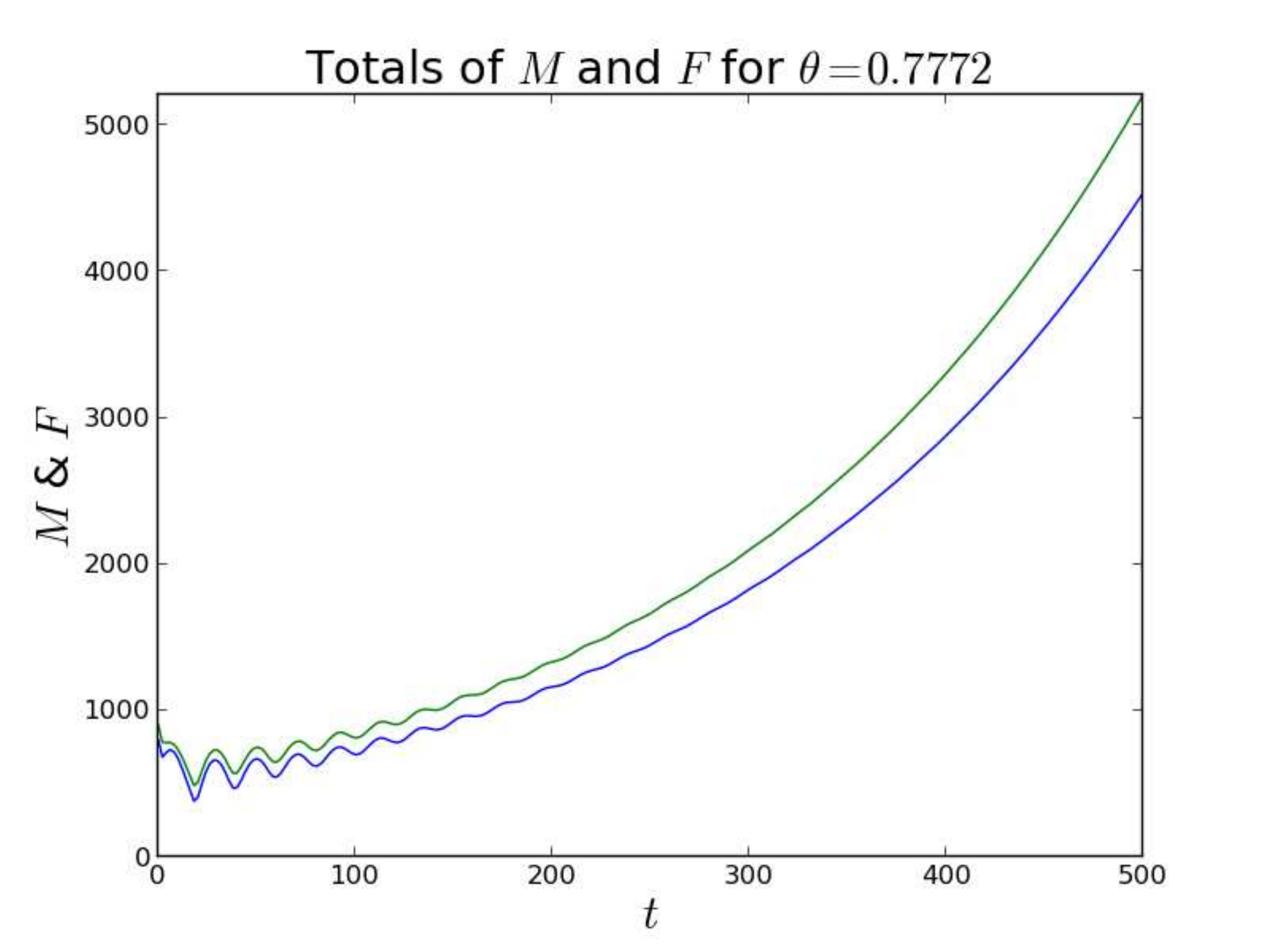}\\
  \caption{Solutions
    to~\eqref{eq:2Bis}--(\ref{eq:parameters-mating})--(\ref{eq:idMF}).
    The integrals $\int_0^{80} M (t,a) \d{a}$ (blue) and $\int_0^{90}
    F (t,a) \d{a}$ (green) as a function of time. Left, for $\theta =
    0.4647$, the population vanishes.  Middle, with $\theta = 0.7067$,
    there is an equilibrium and right, for $\theta = 0.7772$, the
    population grows exponentially. The latter value gives the maximal
    average natality rate, see Figure~\ref{fig:CostMF}. }
  \label{fig:contourMF}
\end{figure}

For completeness, we precise that the numerical integration above was
obtained using a Lax--Friedrichs algorithm,
see~\cite[\S~12.5]{LeVequeBook2002}, with space mesh $\Delta a =
0.04167$.

\subsection{Management of a Biological Resource}
\label{subs:Bio}

In biological resource management, one typically rears/breeds a
species up to a suitable stage, then part of the population is sold
and part is used for reproduction. The
equations~\eqref{eq:1}--\eqref{eq:B} comprehend this
situation. Indeed, call $J = J (t,a)$ the density of the juveniles at
time $t$ of age or size $a$. Juveniles reaching the age/size $\bar a$
are then selected. The density $S = S (t,a)$ refers to those
individuals that are going to be sold, while $R = R (t,a)$ stands for
the density of those reserved for reproduction purposes. One is thus
lead to the following model, defined on the structure in
Figure~\ref{fig:stru2}, right:
\begin{equation}
  \label{eq:6}
  \left\{
    \begin{array}{l@{\qquad}r@{\,}c@{\,}l}
      \partial_t J + \partial_a\left(g_J (t,a) \, J\right)
      =
      d_J (t,a) \, J
      &
      (t,a) & \in & \reali^+ \times [0, \bar a]
      \\
      \partial_t S + \partial_a\left(g_S (t,a) \, S\right)
      =
      d_S (t,a) \, S
      &
      (t,a) & \in & \reali^+ \times \left[\bar a, +\infty\right[
      \\
      \partial_t R + \partial_a\left(g_R (t,a) \, R\right)
      =
      d_R (t,a) \, R
      &
      (t,a) & \in & \reali^+ \times \left[\bar a, +\infty\right[
      \\
      g_J (t, 0)\, J (t, 0)
      =
      \beta\left(\int_{\bar a}^{a_{\max}} R (t,x) \d{x} \right)
      &
      t & \in & \reali^+
      \\
      g_S (t, \bar a) \, S (t, \bar a)
      =
      \eta \, g_J (t, \bar a) \, J (t,\bar a)
      &
      t & \in & \reali^+
      \\
      g_R (t, \bar a) \, R (t, \bar a)
      =
      (1-\eta) \,  g_J (t, \bar a) \, J (t,\bar a)
      &
      t & \in & \reali^+
      \\
      J (0, a) = J_o (a)
      &
      a & \in & [0,\bar a]
      \\
      S (0, a) = S_o (a)
      &
      a & \in & \left[\bar a, +\infty\right[
      \\
      R (0, a) = R_o (a)
      &
      a & \in & \left[\bar a, +\infty\right[ \,.
    \end{array}
  \right.
\end{equation}
Above, we used the obvious notation for the growth and mortality
functions $g_J, g_S, g_R$ and $d_J, d_S, d_R$. The birth rate is
described through the function $\beta$. A key role is played by the
parameter $\eta \in [0,1]$ which quantifies the percentage of
juveniles selected for the market.

System~\eqref{eq:6} fits into~\eqref{eq:1}--\eqref{eq:B} setting
\begin{displaymath}
  \begin{array}{@{}r@{\;}c@{\;}l@{\quad} r@{\;}c@{\;}l@{\quad} r@{\;}c@{\;}l@{}}
    u_1 (t,s)& = & J (t,x)
    &
    g_1 (t,x)& = & g_J (t,x)
    &
    \alpha_1 (t,w_1,w_2,w_3) & = & 0
    \\
    u_2 (t,x)& = & S (t, x+\bar a)
    &
    g_2 (t,x)& = & g_S (t, x+\bar a)
    &
    \alpha_2 (t,w_1,w_2,w_3) & = & \eta \, w_1 \, g_1 (t,\bar x_1)
    \\
    u_3 (t,x)& = & R (t,x+\bar a)
    &
    g_3 (t,x)& = & g_R (t,x+\bar a)
    &
    \alpha_3 (t,w_1,w_2,w_3) & = & (1-\eta) \, w_1 \, g_1 (t,\bar x_1)
    \\
    &  &
    &
    d_1 (t,x)& = & d_J (t,x)
    &
    \beta_1 (w_1,w_2,w_3) & = & \beta(w_3)
    \\
    \bar x_1 & = & \bar a
    &
    d_2 (t,x)& = & d_S (t,x+\bar a)
    &
    \beta_2 (w_1,w_2,w_3) & = & 0
    \\
    I_3 & = & [\bar a, a_{\max}]
    &
    d_3 (t,x)& = & d_R (t, x+\bar a)
    &
    \beta_3 (w_1,w_2,w_3) & = & 0
  \end{array}
\end{displaymath}

\begin{corollary}
  \label{cor:JRS}
  Let $g_J,g_S,g_R$ satisfy~\textbf{(g)} for suitable $\check g, \hat
  g \in \reali^+$ with $\hat g > \check g > 0$.  Let $d_J,d_S,d_R$
  satisfy~\textbf{(d)}. Let $\beta \in \C{0,1} (\reali^+; \reali)$ be
  such that $\beta (0) = 0$. For any $\eta \in [0,1]$ and any initial
  data $J_o \in \BV ([0,\bar a]; \reali^+)$ and $S_o, R_o \in (\L1
  \cap \BV) (\left[\bar a, +\infty\right[, \reali^+)$,
  system~\eqref{eq:6} admits a unique non negative solution and the
  stability estimates in Theorem~\ref{thm:wp} apply.
\end{corollary}

A natural question based on model~\eqref{eq:6} is: find the optimal
percentage $\eta$ of juveniles that have to be chosen for the
market. To this aim, we postulate simple, though reasonable, cost and
gain functionals
\begin{equation}
  \label{eq:CostJSR}
  \!\!\!
  \begin{array}{@{}rcl@{}}
    \mathcal{C} (\eta;T)
    & = &
    \displaystyle
    \int_0^T
    \left[
      \int_0^{\bar a} C_J (a) \, J (t,a) \d{a}
      +
      \int_{\bar a}^{a_{\max}}
      \left[
        C_S (a) \, S (t,a)
        +
        C_R (a) \, R (t,a)
      \right]
      \d{a}
    \right]
    \d{t} ,
    \\
    \mathcal{G} (\eta;T)
    & = &
    \displaystyle
    \int_0^T
    \int_{\bar a}^{a_{\max}}
    G (a) \, S (t,a) \d{a} \, \d{t}\,.
  \end{array}
\end{equation}
Here, $C_J (a)$ is the unit cost to grow a juvenile at age $a$, and
similarly $C_R$, $C_S$ are the costs for the other two groups. The
gain obtained selling an adult at age $a$ is $G (a)$. We denoted by $J
= J (t,a)$, $S = S (t,a)$ and $R = R (t,a)$ the solution
to~\eqref{eq:6} with initial datum $J_o$ and $S_o=0$, $R_o=0$ and with
the selection parameter $\eta$. A direct consequence of
Theorem~\ref{thm:Stab} is the following corollary.

\begin{corollary}
  \label{cor:JRS2}
  In the same assumptions of Corollary~\ref{cor:JRS}, for any $T \in
  \pint{\reali}^+$ and any $J_o, S_o, R_o \in (\L1 \cap \BV)
  (\reali^+; \reali)$, there exists an optimal choice $\eta_*$ such
  that
  \begin{displaymath}
    \mathcal{G} (\eta_*;T) - \mathcal{C} (\eta_*;T)
    =
    \max_{\eta \in [0,1]}
    \left(
      \mathcal{G} (\eta;T) - \mathcal{C} (\eta;T)
    \right) \,.
  \end{displaymath}
\end{corollary}

\noindent The proof relies on {Weierstra\ss} Theorem, exactly as that
of Corollary~\ref{cor:mating2} and is here omitted.

\begin{figure}[h!]
  \centering
  \includegraphics[width=0.5\textwidth]{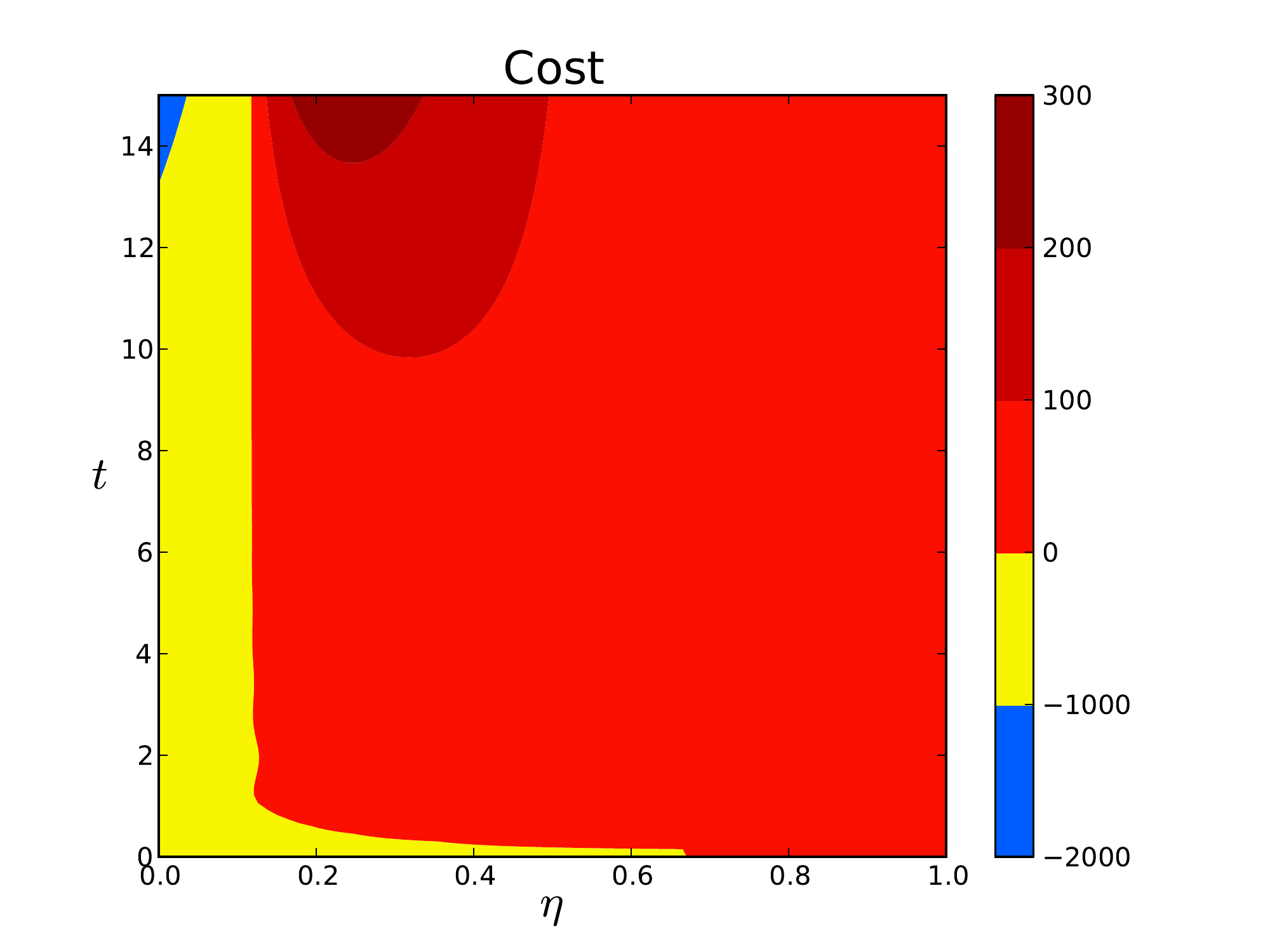}%
  \includegraphics[width=0.5\textwidth]{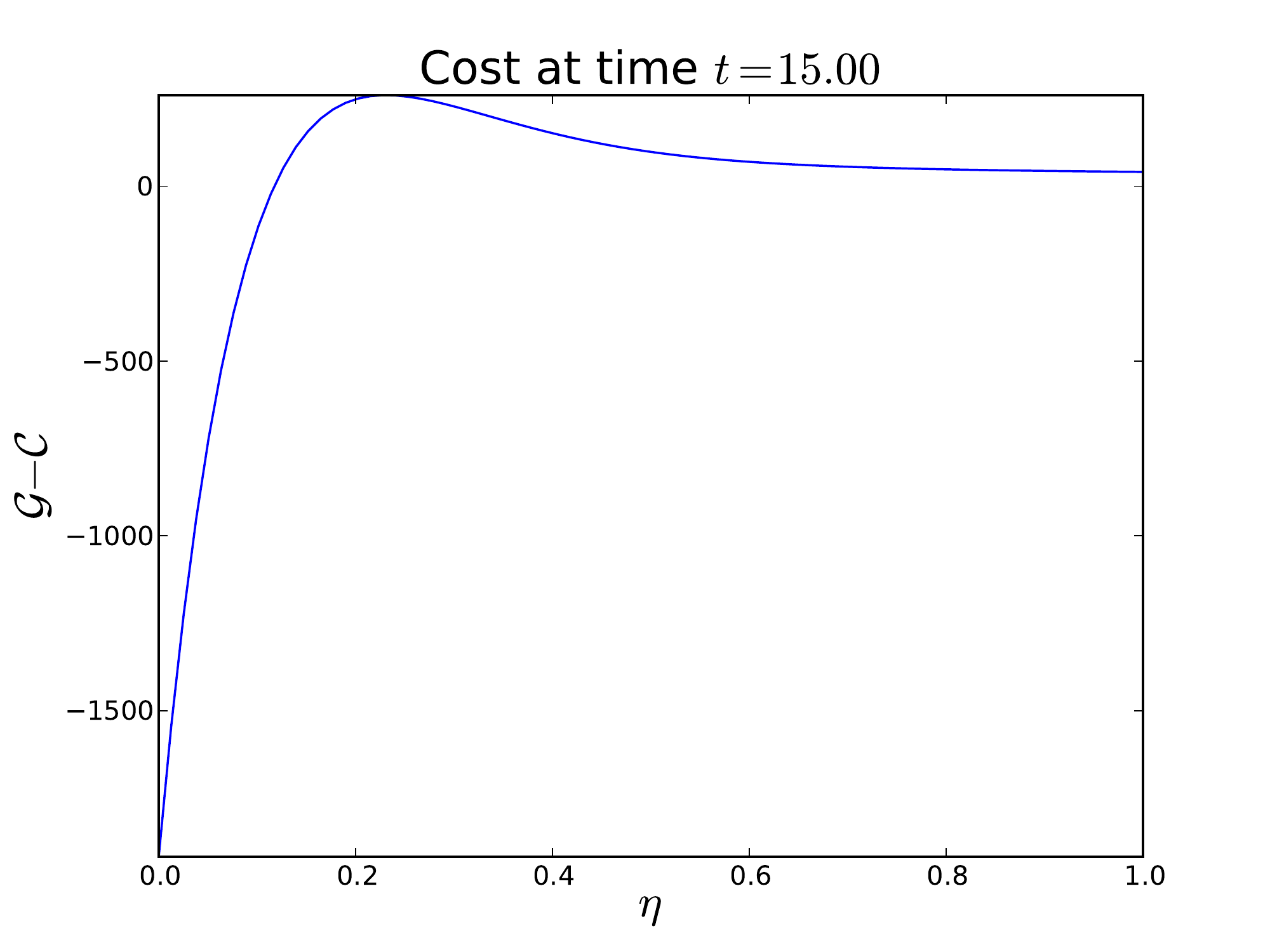}
  \caption{Cost~\eqref{eq:CostJSR} along the solutions
    to~\eqref{eq:6}--\eqref{eq:8} with initial
    datum~\eqref{eq:idJSR}. Left, as a function of $\eta$ (on the
    horizontal axis) and $t$ (on the vertical axis). Right, as a
    function of $\eta$ at time $t=15.00$. Recall that for $\eta = 0$,
    all individual are kept for reproduction and no one is sold. On
    the contrary, for $\eta = 1$, they are all sold and no one is kept
    for reproduction.}
  \label{fig:CostJSR}
\end{figure}

As a specific example, we consider the situation identified by the
following choices of functions and parameters
in~\eqref{eq:6}--\eqref{eq:CostJSR}:
\begin{equation}
  \label{eq:8}
  \begin{array}{@{}rcl@{\qquad}rcl@{\qquad}rcl@{\qquad}rcl@{}}
    g_J (t,a) & = & 1
    &
    d_J (t,a) & = & 0
    &
    C_J (t,a) & = & a
    &
    \beta (w) & = & 2\, w
    \\
    g_S (t,a) & = & 1
    &
    d_S (t,a) & = & -\frac{a-\bar a}{2}
    &
    C_S (t,a) & = & 0
    &
    G (t,a) & = & 10
    \\
    g_R (t,a) & = & 1
    &
    d_R (t,a) & = & -\frac{a-\bar a}{2}
    &
    C_R (t,a) & = & 0.5
    &
    [\bar a, a_{\max}] & = & [1, 2]
  \end{array}
\end{equation}
and we consider $\eta$ as a control parameter in $[0,1]$.  As initial
datum we choose
\begin{equation}
  \label{eq:idJSR}
  J_o (a) = 5
  \,,\quad
  S_o (a) = 0
  \,,\quad
  R_o (a) = 0 \,.
\end{equation}
The graph of the cost $\mathcal G(\eta; T) - \mathcal C(\eta; T)$
(see~\eqref{eq:CostJSR}) for $T = 15$ with respect to $\lambda$ is in
Figure~\ref{fig:CostJSR}.
The outcome shows a reasonable qualitative behavior. As $\eta \to 0$,
nothing is sold, all population members are kept for reproduction, the
population increases exponentially as also does the functional
$\mathcal G(\eta;T) - \mathcal C(\eta;T)$. On the contrary, for $\eta
\to 1$, all population members are immediately sold giving a positive
gain and the population vanishes as also $\mathcal G -\mathcal
C$. Near to $\eta \approx 0.23$ there is an optimal balance, given the
chosen unitary costs and gain~\eqref{eq:CostJSR}--\eqref{eq:8}.

\begin{figure}[h!]
  \centering
  \includegraphics[width=0.33\textwidth,trim=20 5 20 5]{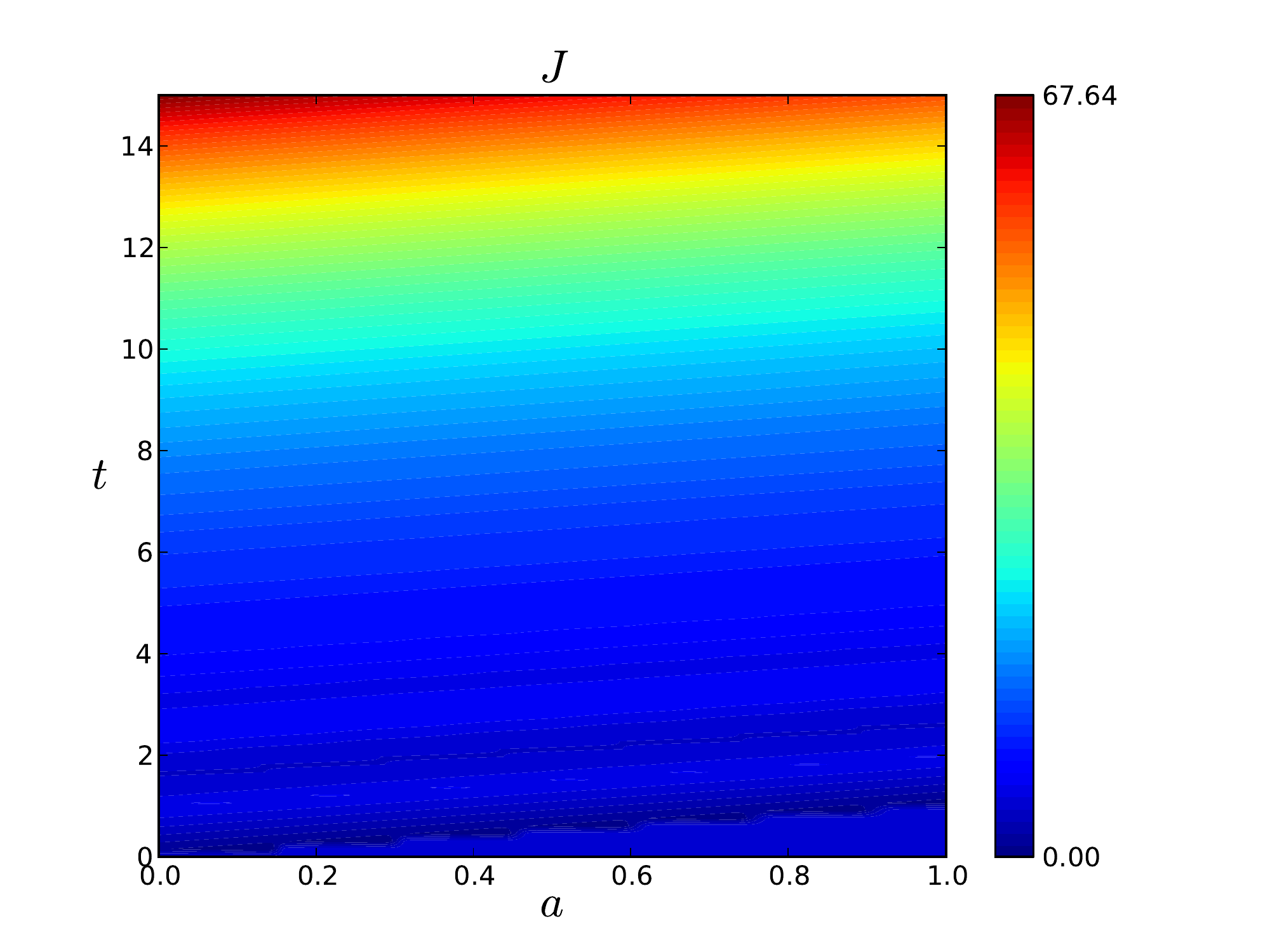}%
  \includegraphics[width=0.33\textwidth,trim=20 5 20 5]{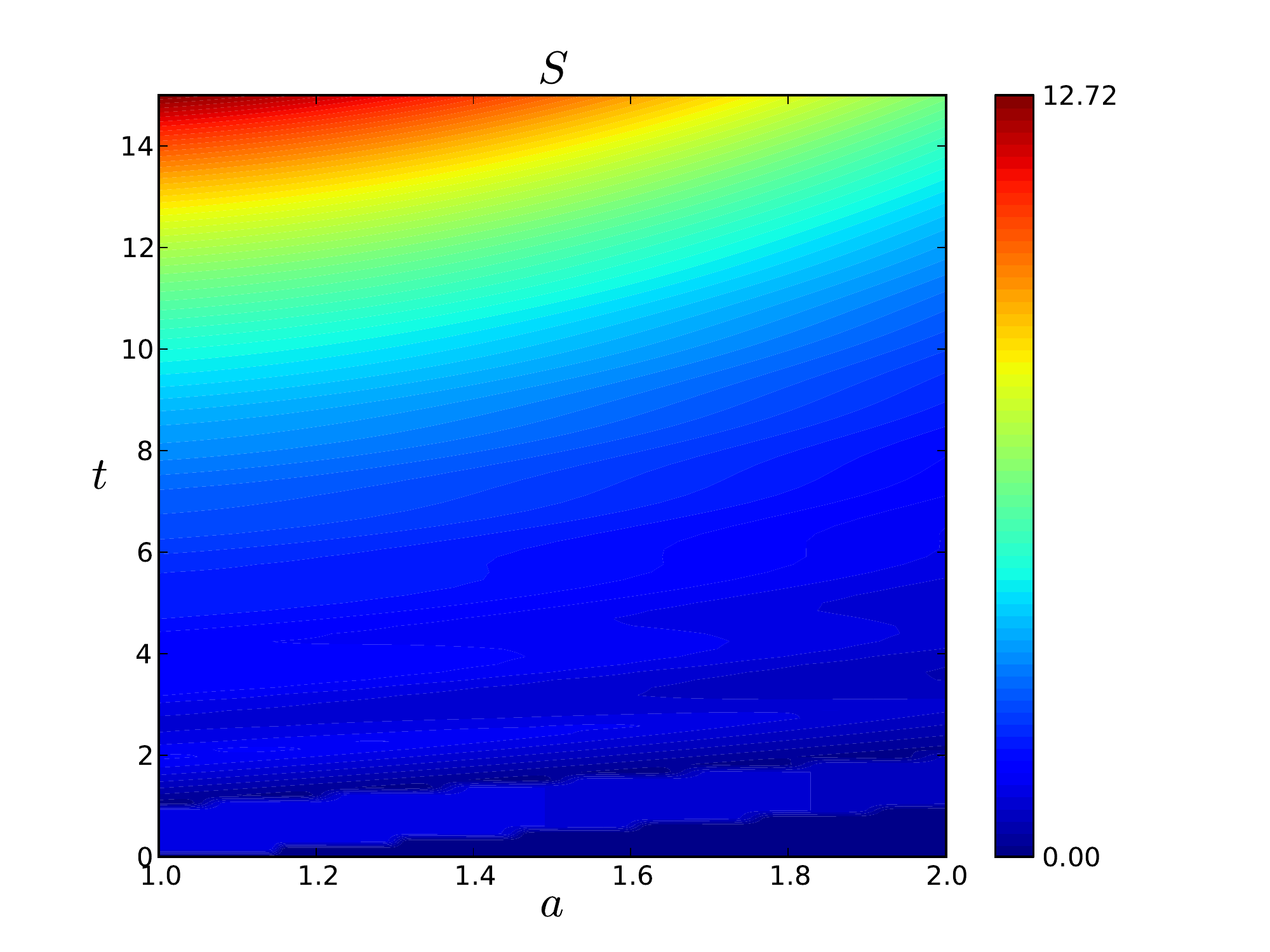}%
  \includegraphics[width=0.33\textwidth,trim=20 5 20 5]{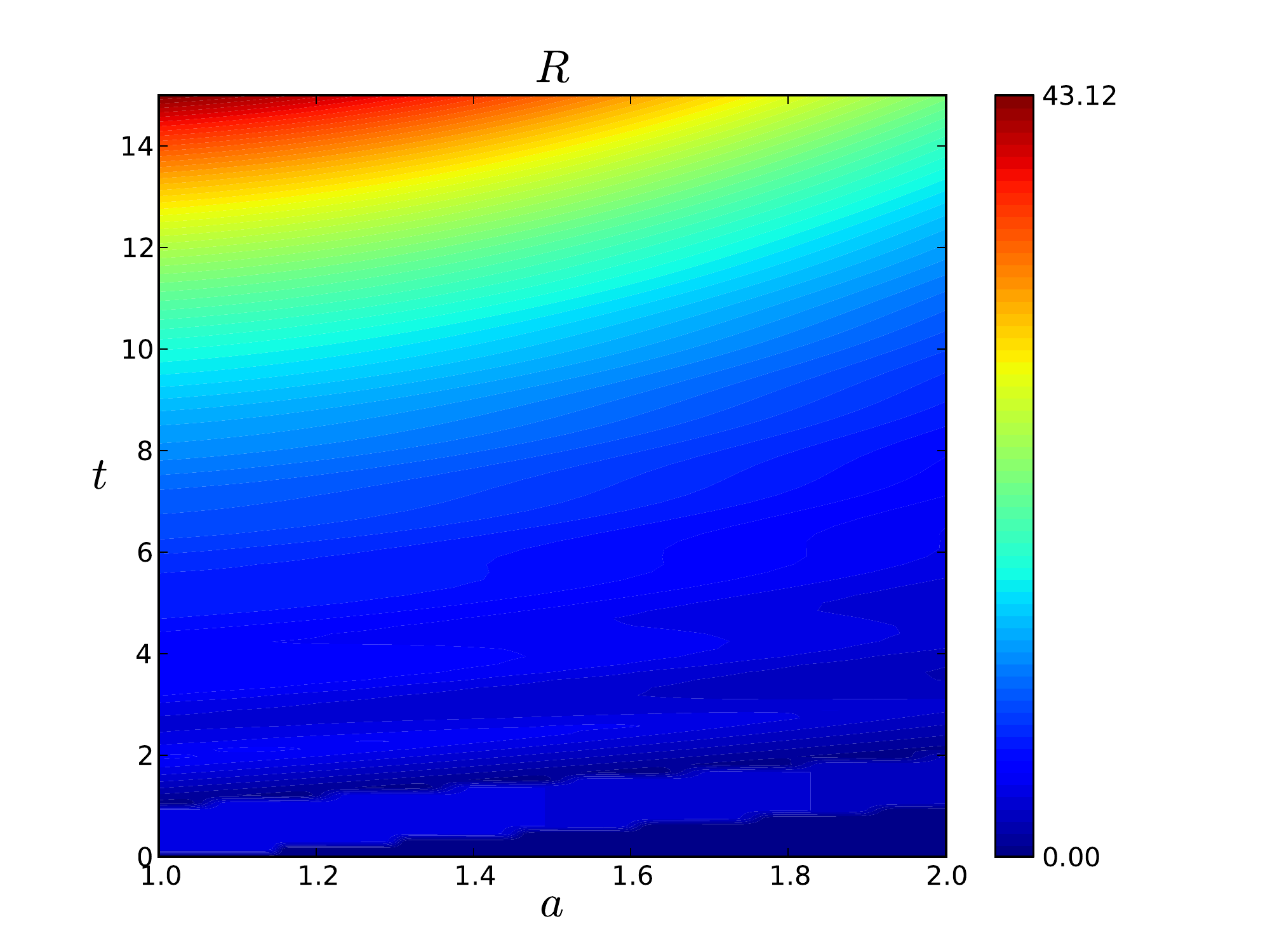}\\
  \includegraphics[width=0.33\textwidth,trim=20 5 20 5]{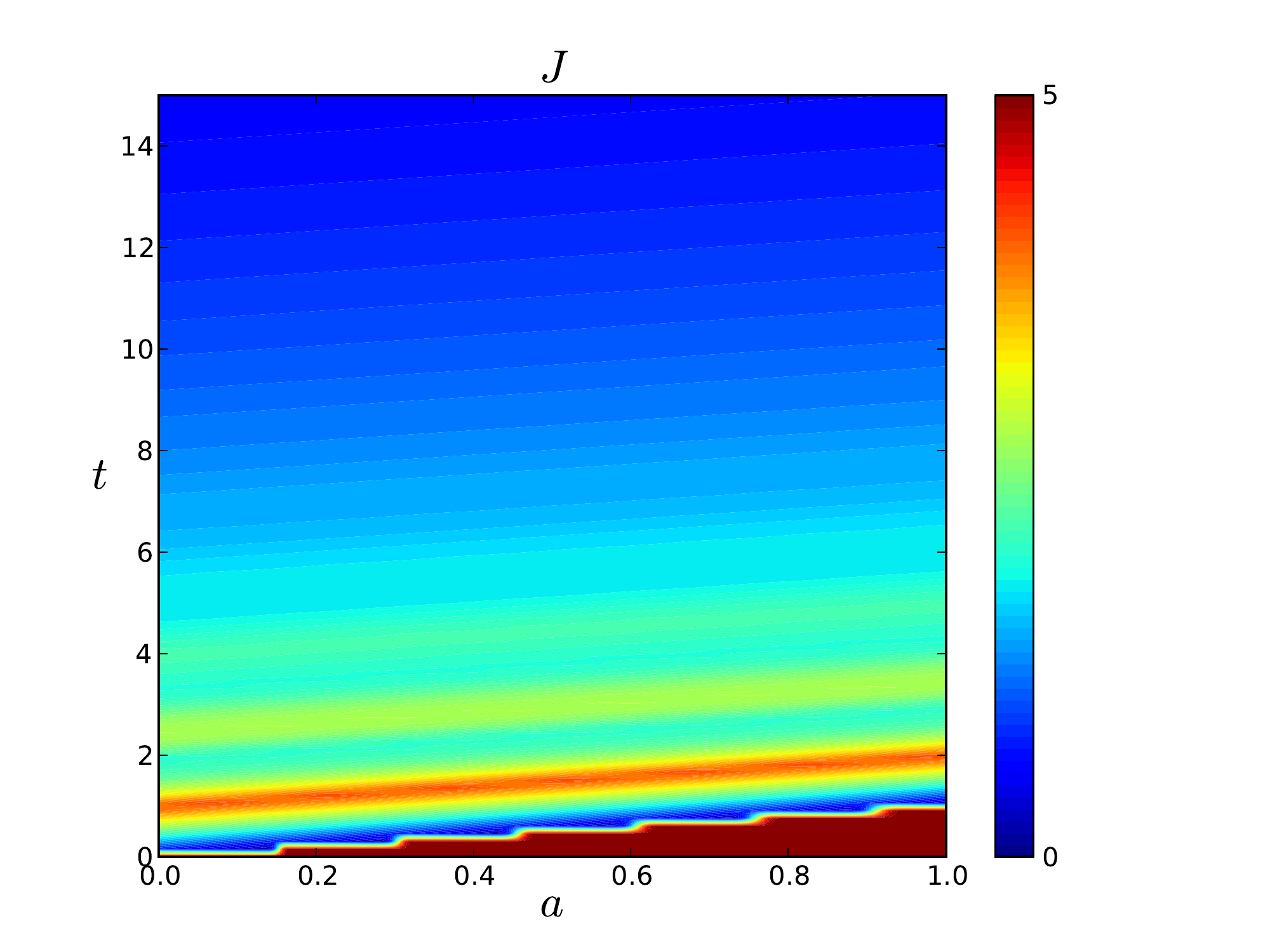}%
  \includegraphics[width=0.33\textwidth,trim=20 5 20 5]{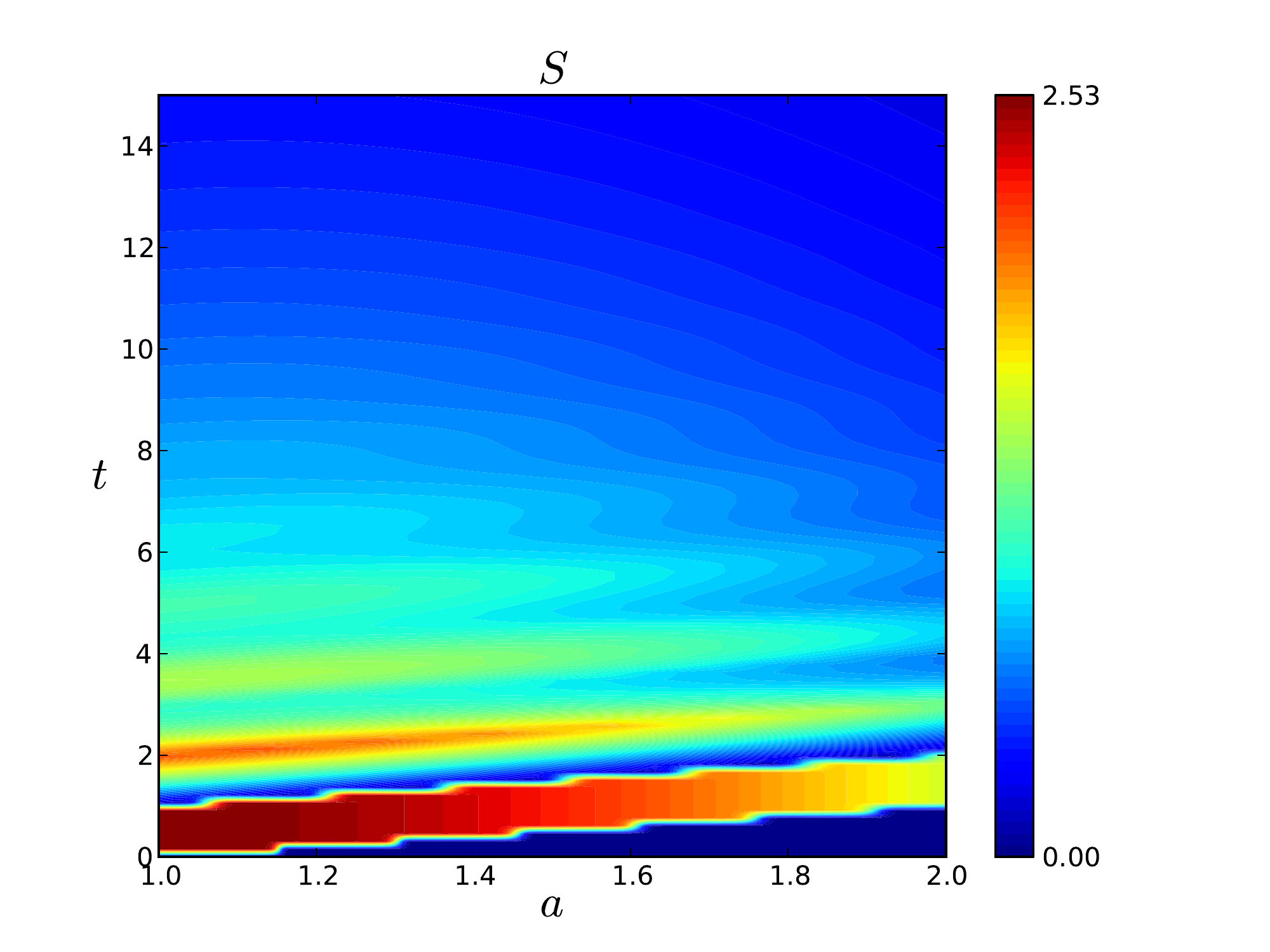}%
  \includegraphics[width=0.33\textwidth,trim=20 5 20 5]{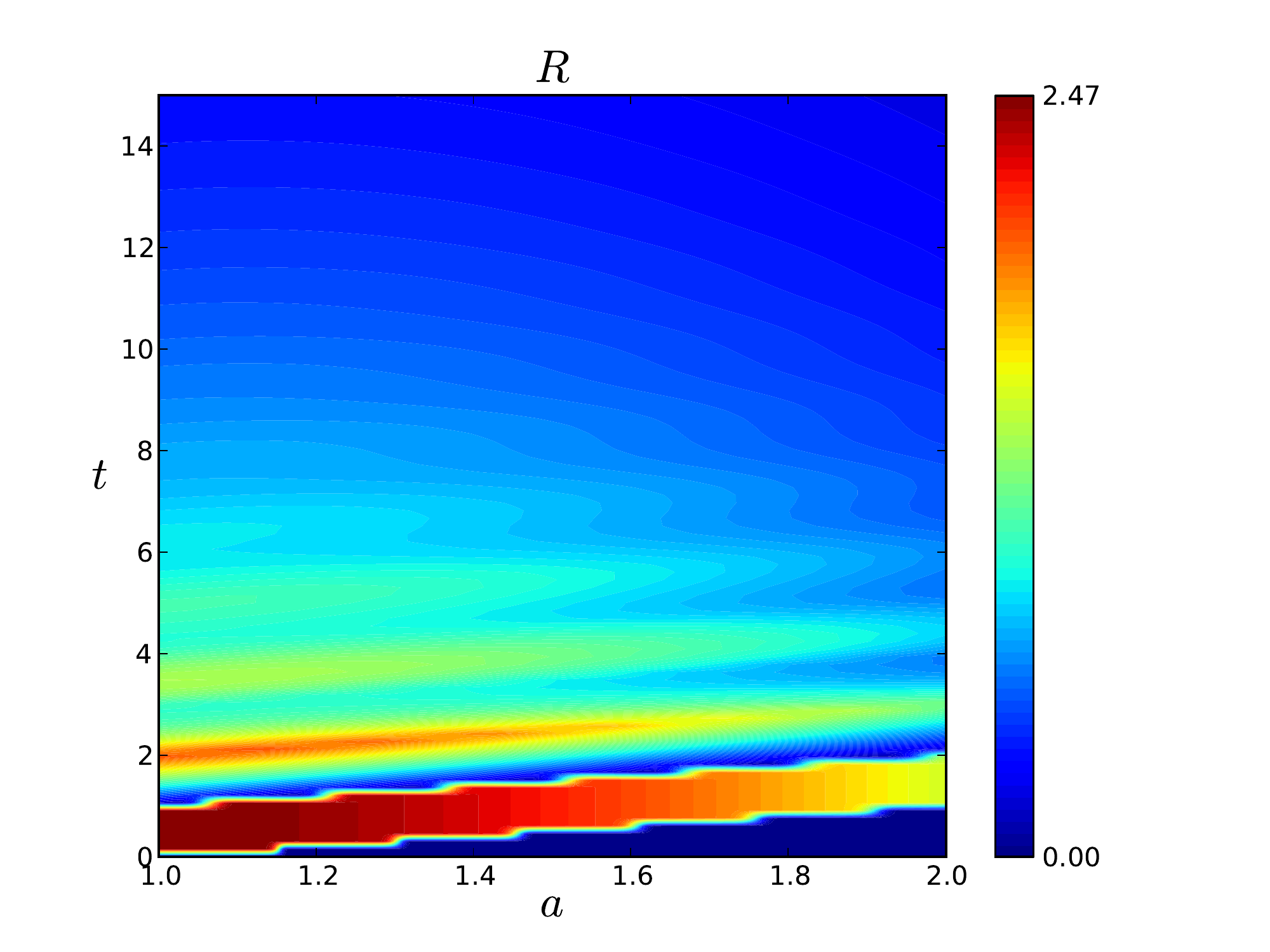}\\
  \includegraphics[width=0.33\textwidth,trim=20 5 20 5]{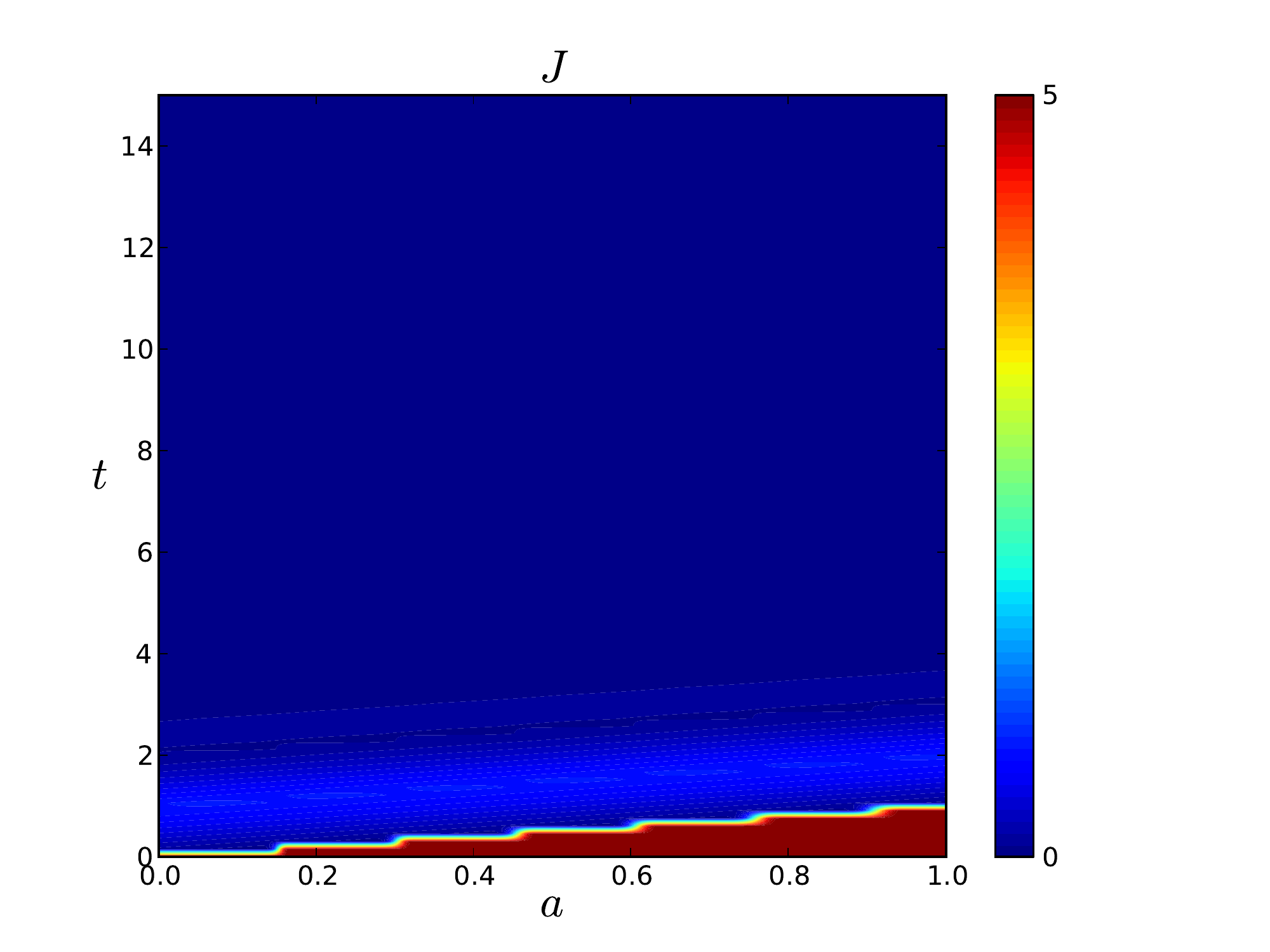}%
  \includegraphics[width=0.33\textwidth,trim=20 5 20 5]{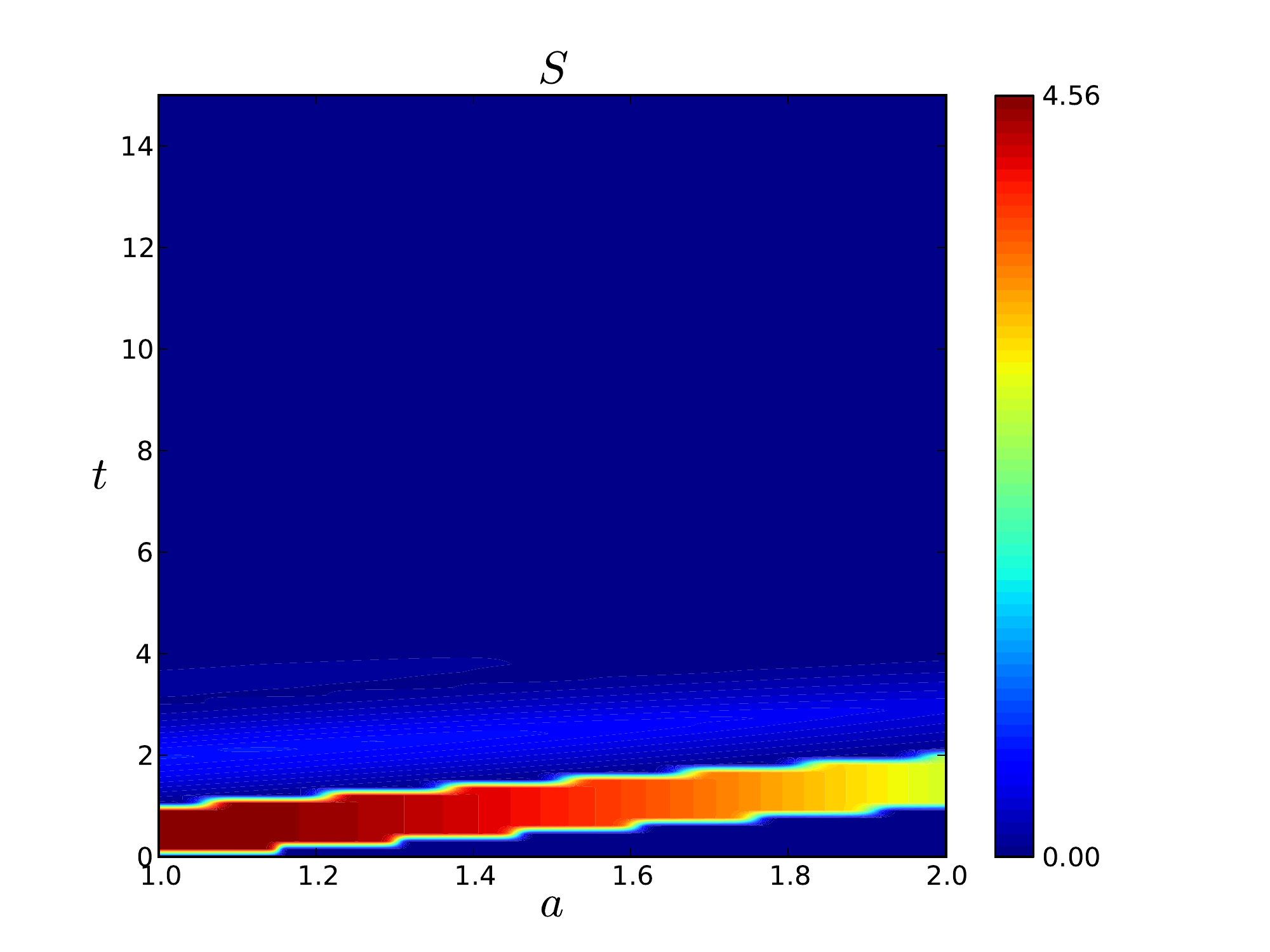}%
  \includegraphics[width=0.33\textwidth,trim=20 5 20 5]{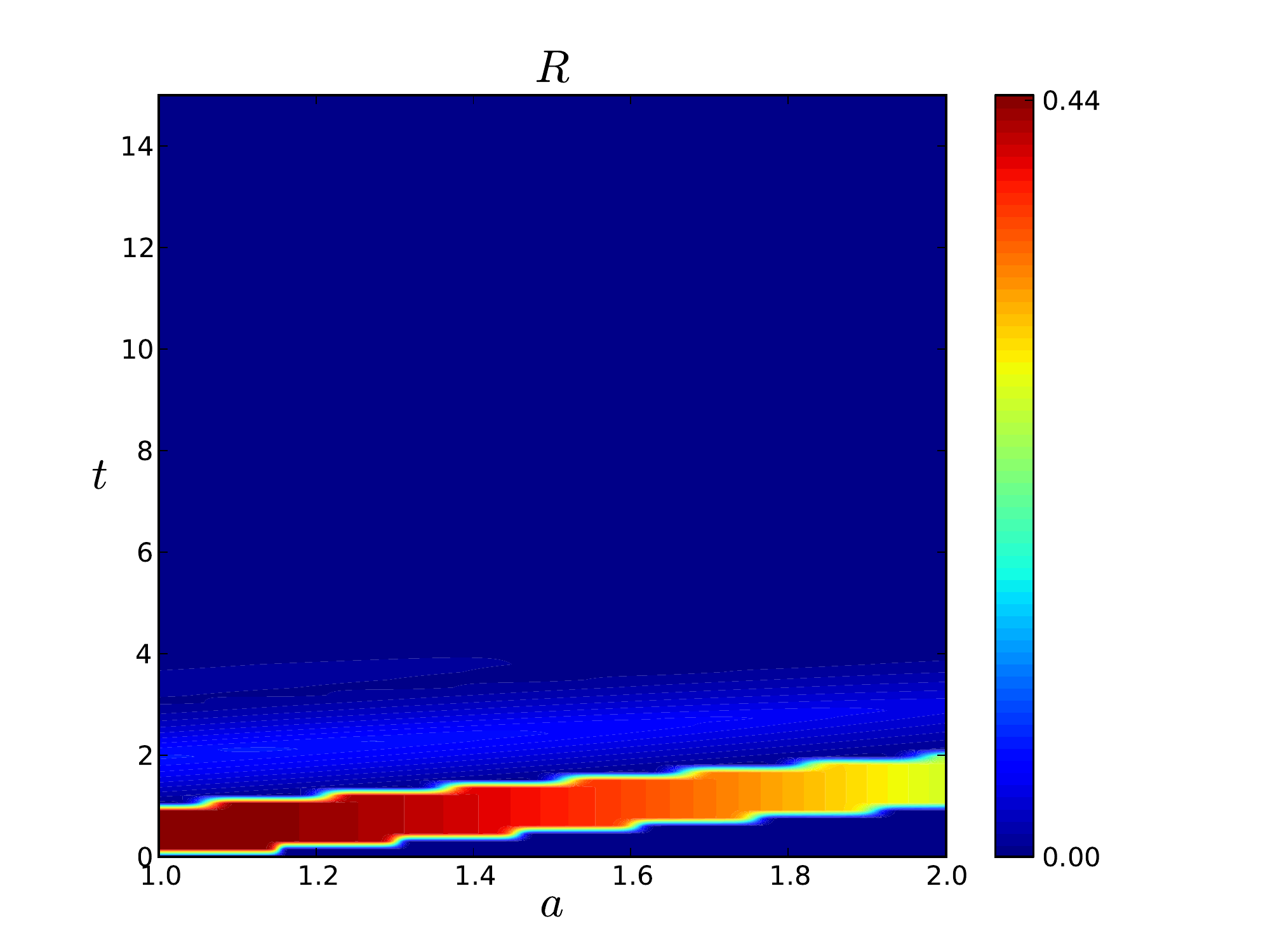}\\
  \caption{Solutions
    to~\eqref{eq:6}--\eqref{eq:8}--\eqref{eq:idJSR}. Time varies along
    the vertical axis and $a$ along the horizontal one. Above, $\eta =
    0.23$ is near to the optimal choice. Middle, $\eta = 0.50$ and,
    below, $\eta = 0.91$.}
  \label{fig:contourJSR}
\end{figure}

With the chosen parameters, the optimal choice for $\eta$ is $\eta_*
\approx 0.23$, which yields a gain of about $260.48$ at time
$t=15$. As expected, different choices of $\eta$ have deep influences
on the solutions to~\eqref{eq:6}--\eqref{eq:8}--\eqref{eq:idJSR}, as
shown in Figure~\ref{fig:contourJSR}.

For completeness, we precise that the numerical integration above was
obtained using a Lax--Friedrichs algorithm,
see~\cite[\S~12.5]{LeVequeBook2002}, with space mesh $\Delta a =
0.001$.

\section{Technical Details}
\label{sec:TD}

Throughout, when $\BV$ functions are considered, we refer to a right
continuous representative. We now recall the following elementary
estimates on $\BV$ functions.
\begin{eqnarray}
  \label{eq:TV1}
  \left.
    \begin{array}{r@{\,}c@{\,}l@{}}
      u & \in & \BV (\reali^+; \reali)
      \\
      w & \in & \BV (\reali^+; \reali)
    \end{array}
  \right\}
  &\Rightarrow &
  \tv (u \, w)
  \leq
  \norma{u}_{\L\infty (\reali^+; \reali)} \tv (w)
  +
  \tv (u) \, \norma{w}_{\L\infty (\reali^+; \reali)}
  \\
  \label{eq:TV2}
  \left.
    \begin{array}{r@{\,}c@{\,}l@{}}
      f & \in & \C{0,1} (\reali; \reali)
      \\
      u & \in & \BV (\reali^+; \reali)
    \end{array}
  \right\}
  &\Rightarrow &
  \tv (f\circ u) \leq \Lip (f) \, \tv (u)
  \\
  \label{eq:TV4}
  \left.
    \begin{array}{@{}r@{\,}c@{\,}l@{}}
      u & \in & \BV (\reali^+; \reali)
      \\
      f & \in & \BV (\reali^+; [\check f +\infty[\,)
      \\
      \check f & > & 0
    \end{array}
  \right\}
  &\Rightarrow &
  \tv\left(\frac{u}{f}\right)
  \leq
  \frac{1}{\check f} \, \tv (u)
  +
  \frac{1}{\check f^2} \, \tv (f) \, \norma{u}_{\L\infty (\reali^+; \reali)}
  \\
  \label{eq:TV3}
  \!\!\!u \in \L1 \!\left(\reali^+; \BV (\reali^+;\reali)\right)
  &\Rightarrow &
  \tv\left(\int_0^t u (\tau,\cdot) \, \d\tau\right)
  \leq
  \int_0^t \tv\!\left(u (\tau)\right) \, \d\tau
  \\
  \label{eq:3}
  \left.
    \begin{array}{r@{\,}c@{\,}l@{}}
      u & \in & \BV (\reali^+; \reali)
      \\
      h & \in & \L\infty(\reali; \reali^+)
    \end{array}
  \right\}
  &\Rightarrow &
  \int_{\reali^+}
  \modulo{u\left(x + h (x)\right) - u (x)} \d{x}
  \leq
  \tv (u) \, \norma{h}_{\L\infty (\reali^+; \reali)} \,.
\end{eqnarray}
Inequality~\eqref{eq:TV1} follows
from~\cite[Formula~(3.10)]{AmbrosioFuscoPallara}. The definition of
total variation directly implies~\eqref{eq:TV2}, \eqref{eq:TV4}
and~\eqref{eq:TV3}.  For a proof of~\eqref{eq:3} see for
instance~\cite[Lemma~2.3]{BressanLectureNotes}.

\begin{proofof}{Lemma~\ref{lem:stability}}
  To verify that~(\ref{eq:12}) solves~(\ref{eq:4}), a standard
  integration along characteristics is sufficient. The
  bounds~\eqref{eq:L1} and~\eqref{eq:Linfty} are an immediate
  consequence of~\eqref{eq:12}.

  Passing to the estimates on the total variation, introduce
  \begin{equation}
    \label{eq:C}
    C = 2 \, \max
    \left\{
      \begin{array}{ll}
        \norma{\partial_x g}_{\L\infty (\reali^+ \times \reali^+; \reali)},
        &
        \norma{\partial_t g}_{\L\infty (\reali^+ \times \reali^+; \reali)},
        \\
        \sup_{t \in \reali^+} \tv\left(g (t,\cdot)\right),
        &
        \sup_{t \in \reali^+} \tv\left(\partial_x g (t,\cdot)\right),
        \\
        \norma{d}_{\L\infty (\reali^+ \times \reali^+; \reali)},
        &
        \sup_{t \in \reali^+} \tv\left(d(t,\cdot)\right)
      \end{array}
    \right\}
  \end{equation}
  which is finite by~\textbf{(g)} and~\textbf{(d)}.

  Consider now the total variation estimates. Using~(\ref{eq:TV1}),
  (\ref{eq:TV2}), (\ref{eq:TV4}), (\ref{eq:TV3}), compute:
  \begin{eqnarray*}
    \tv\left(u (t)\right)
    & = &
    \tv\left(u (t, \cdot), [0, \gamma (t)]\right)
    +
    \tv\left(u (t, \cdot), \left[\gamma (t), +\infty\right[\right)
    \\
    & \leq &
    \tv \left(b (\cdot)/g (\cdot , 0); [0,t]\right) e^{C t}
    +
    \frac{\norma{b}_{\L\infty ([0,t]; \reali)}}{\check g} e^{C t}
    +
    \tv (u_o) e^{C t} + \norma{u_o}_{\L\infty (\reali^+; \reali)} e^{C t}
    \\
    & \leq &
    \left(
      \norma{u_o}_{\L\infty (\reali^+; \reali)}
      +
      \tv (u_o)
      +
      \frac{1}{\check g}
      \left(
        \left(1+\frac{C}{\check g}\right) \,
        \norma{b}_{\L\infty ([0,t]; \reali)}
        +
        \tv (b;[0,t])
      \right)
    \right)
    e^{C t}
  \end{eqnarray*}
  proving~\eqref{eq:tvx}.  The bound~\eqref{eq:tvx} directly follows
  from~\eqref{eq:12}, using~(\ref{eq:TV1}), (\ref{eq:TV2})
  and~(\ref{eq:TV4}).  We exploit now~(\ref{eq:TV3})
  and~\cite[Definition~3.4]{AmbrosioFuscoPallara}, and in the lines
  below, for typographical reasons, we denote by $\uno$ the real
  interval $[-1,1]$.
  \begin{eqnarray*}
    \!\!\! & &
    \tv\left(\int_I u (\cdot, x) \d{x}; [0,t]\right)
    \\
    \!\!\! & = &
    \sup
    \left\{
      \int_0^t \int_I
      u (\tau,x) \, \d{x} \, \partial_t \phi (\tau) \, \d\tau
      \colon \phi \in \Cc1 (\left]0,t\right[;\uno)
    \right\}
    \\
    \!\!\! & \leq &
    \sup
    \left\{
      \int_0^t \int_I
      u (\tau,x) \, \psi (x) \, \partial_t \phi (\tau) \, \d{x} \, \d\tau
      \colon
      \begin{array}{@{}r@{\,}c@{\,}l@{}}
        \phi & \in & \Cc1 (\left]0,t\right[;\uno)
        \\
        \psi & \in & \Cc1 (\pint{I}; \uno)
      \end{array}
    \right\}
    \\
    \!\!\! & = &
    \sup
    \left\{
      \int_0^t \int_I
      \left(
        g (\tau,x) u (\tau,x) \partial_x \psi (x) \phi (\tau)
        +
        d (\tau,x) u (\tau,x) \psi (x) \phi (\tau)
      \right) \d{x} \d\tau
      \colon
      \begin{array}{@{}r@{\,}c@{\,}l@{}}
        \phi & \in & \Cc1 (\left]0,t\right[;\uno)
        \\
        \psi & \in & \Cc1 (\pint{I}; \uno)
      \end{array}
    \right\}
    \\
    \!\!\! & = &
    \sup
    \left\{
      \int_0^t \int_I
      g (\tau,x) \, u (\tau,x) \, \partial_x \psi (x) \, \phi (\tau)
      \, \d{x} \, \d\tau
      \colon
      \begin{array}{@{}r@{\,}c@{\,}l@{}}
        \phi & \in & \Cc1 (\left]0,t\right[;\uno)
        \\
        \psi & \in & \Cc1 (\pint{I}; \uno)
      \end{array}
    \right\}
    \\
    \!\!\! & &
    +
    \sup
    \left\{
      \int_0^t \int_I
      d (\tau,x) \, u (\tau,x) \, \psi (x) \, \phi (\tau)
      \, \d{x} \, \d\tau
      \colon
      \begin{array}{@{}r@{\,}c@{\,}l@{}}
        \phi & \in & \Cc1 (\left]0,t\right[;\uno)
        \\
        \psi & \in & \Cc1 (\pint{I}; \uno)
      \end{array}
    \right\}
    \\
    \!\!\! & \leq &
    \sup
    \left\{
      \int_0^t
      \sup
      \left\{
        \int_I
        g (\tau,x) \, u (\tau,x) \, \partial_x \psi (x)
        \, \d{x}
        \colon
        \psi \in \Cc1 (\pint{I}; \uno)
      \right\}
      \, \phi (\tau)
      \, \d\tau
      \colon
      \phi \in \Cc1 (\left]0,t\right[;\uno)
    \right\}
    \\
    \!\!\! & &
    +
    \sup
    \left\{
      \int_0^t
      \sup
      \left\{
        \int_I
        d (\tau,x) \, u (\tau,x) \, \psi (x)
        \, \d{x}
        \colon
        \psi\in \Cc1 (\pint{I}; \uno)
      \right\}
      \, \phi (\tau)
      \, \d\tau
      \colon
      \phi \in \Cc1 (\left]0,t\right[;\uno)
    \right\}
    \\
    \!\!\! & = &
    \sup
    \left\{
      \int_0^t
      \tv\left( g (\tau,\cdot) \, u (\tau,\cdot ) \right)
      \, \phi (\tau)
      \, \d\tau
      \colon
      \phi \in \Cc1 (\left]0,t\right[;\uno)
    \right\}
    \\
    \!\!\! & &
    +
    \sup
    \left\{
      \int_0^t
      \tv\left(d (\tau,\cdot) \, u (\tau,\cdot)\right)
      \, \phi (\tau)
      \, \d\tau
      \colon
      \phi \in \Cc1 (\left]0,t\right[;\uno)
    \right\}
    \\
    \!\!\! & \leq &
    \int_0^t
    \left(
      \tv\left( g (\tau,\cdot) \, u (\tau,\cdot ) \right)
      +
      \tv\left(d (\tau,\cdot) \, u (\tau,\cdot)\right)
    \right)
    \, \d\tau
  \end{eqnarray*}
  Apply now~\eqref{eq:TV1} to obtain:
  \begin{eqnarray*}
    \tv\left(\int_I u (\cdot, x) \d{x}; [0,t]\right)
    & \leq &
    \int_0^t
    \left(
      \tv\left(g (\tau, \cdot)\right)
      +
      \tv\left(d (\tau, \cdot)\right)
    \right)
    \norma{u (\tau)}_{\L\infty (\reali^+; \reali)}
    \d\tau
    \\
    & &
    +
    \int_0^t
    \left(
      \norma{g (\tau)}_{\L\infty (\reali^+; \reali)}
      +
      \norma{d (\tau)}_{\L\infty (\reali^+; \reali)}
    \right)
    \tv\left(u (\tau)\right)
    \d\tau
    \\
    & \leq &
    2 \, C \int_0^t
    \left(
      \norma{u (\tau)}_{\L\infty (\reali^+; \reali)}
      +
      \tv\left(u (\tau)\right)
    \right)
    \d\tau \,,
  \end{eqnarray*}
  completing the proof of~\eqref{eq:tvi}. Concerning the stability
  bounds, \eqref{eq:12} implies
  \begin{eqnarray*}
    \int_0^{\gamma (t)} \modulo{u (t,x)} \d{x}
    & \leq &
    \frac{1}{\check g} \int_0^t \modulo{b (\tau)} \d\tau
    +
    \int_0^t \int_{0}^{\gamma (t)}
    \modulo{d (\tau,x) \, u (\tau, x)} \d{x} \d\tau
    \\
    \int_{\gamma (t)}^{+\infty} \modulo{u (t,x)} \d{x}
    & \leq &
    \int_{0}^{+\infty} \modulo{u_o (x)} \d{x}
    +
    \int_0^t \int_{\gamma (t)}^{+\infty}
    \modulo{d (\tau,x) \, u (\tau, x)}  \d{x} \d\tau
    \\
    \norma{u (t)}_{\L1 (\reali^+; \reali)}
    &\leq &
    \norma{u_o}_{\L1 (\reali^+; \reali)}
    +
    \frac{1}{\check g} \, \norma{b}_{\L1 ([0,t]; \reali)}
    +
    C \int_0^t \norma{u (\tau)}_{\L1 (\reali^+; \reali)} \d\tau \,.
  \end{eqnarray*}
  An application of Gronwall Lemma yields the desired
  estimate~\eqref{eq:uffa}. Finally, the monotonicity
  property~\eqref{eq:mono} directly follows from~\eqref{eq:12}.

  To prove~\eqref{eq:lip-dependence}, fix $t',t'' \in \reali^+$ with
  $t' < t''$. Then,
  \begin{eqnarray}
    \label{eq:LipT1}
    \norma{u (t'') - u (t')}_{\L1 (\reali^+;\reali)}
    & = &
    \int_0^{\gamma (t')} \modulo{u (t'',x) - u (t',x)} \d{x}
    \\
    \label{eq:LipT2}
    & &
    +
    \int_{\gamma (t')}^{\gamma (t'')} \modulo{u (t'',x) - u (t',x)} \d{x}
    \\
    \label{eq:LipT3}
    & &
    +
    \int_{\gamma (t'')}^{+\infty} \modulo{u (t'',x) - u (t',x)} \d{x}
  \end{eqnarray}
  and we deal with the three terms separately, using~\eqref{eq:12} as
  follows. Begin with~\eqref{eq:LipT1}:
  \begin{eqnarray*}
    & &
    \int_0^{\gamma (t')} \modulo{u (t',x) - u (t'',x)} \d{x}
    \\
    & \leq &
    \int_0^{\gamma (t')}
    \Bigg|
    \frac{b\left(T (0;t',x)\right)}{g\left(T (0;t',x),0\right)}
    \exp\left(
      \int_{T (0;t',x)}^{t'}
      \left(
        d \left(\tau, X (\tau; t', x)\right)
        -
        \partial_x g \left(\tau, X (\tau; t', x)\right)
      \right)
      \d{\tau}
    \right)
    \\
    & &
    \qquad
    -
    \frac{b\left(T (0;t'',x)\right)}{g\left(T (0;t'',x),0\right)}
    \exp\left(
      \int_{T (0;t'',x)}^{t''}
      \left(
        d \left(\tau, X (\tau; t'', x)\right)
        -
        \partial_x g \left(\tau, X (\tau; t'', x)\right)
      \right)
      \d{\tau}
    \right)
    \Bigg|
    \d{x}
    \\
    & \leq &
    \int_0^{\gamma (t')}
    \modulo{
      \frac{b\left(T (0;t',x)\right)}{g\left(T (0;t',x),0\right)}
      -
      \frac{b\left(T (0;t'',x)\right)}{g\left(T (0;t'',x),0\right)}
    }
    \\
    & &
    \qquad
    \times
    \exp\left(
      \int_{T (0;t',x)}^{t'}
      \left(
        d \left(\tau, X (\tau; t', x)\right)
        -
        \partial_x g \left(\tau, X (\tau; t', x)\right)
      \right)
      \d{\tau}
    \right)
    \d{x}
    \\
    & &
    +
    \int_0^{\gamma (t')}
    \modulo{\frac{b\left(T (0;t'',x)\right)}{g\left(T (0;t'',x),0\right)}}
    \\
    & &
    \qquad
    \times
    \Bigg|
    \exp\left(
      \int_{T (0;t',x)}^{t'}
      \left(
        d \left(\tau, X (\tau; t', x)\right)
        -
        \partial_x g \left(\tau, X (\tau; t', x)\right)
      \right)
      \d{\tau}
    \right)
    \\
    & &
    \qquad\qquad
    -
    \exp\left(
      \int_{T (0;t'',x)}^{t''}
      \left(
        d \left(\tau, X (\tau; t'', x)\right)
        -
        \partial_x g \left(\tau, X (\tau; t'', x)\right)
      \right)
      \d{\tau}
    \right)
    \Bigg|
    \d{x}
    \\
    & \leq &
    \hat g \, t' \,
    \left(
      \frac{\tv (b; [0,t''])}{\check g}
      +
      \norma{b}_{\L\infty ([0,t'']; \reali)}
      \frac{\tv\left(g (\cdot, 0); [0,t'']\right)}{\check g^2}
    \right)
    \, (t''-t') \, e^{2C\, t'}
    \\
    & &
    +
    \hat g \, t'  \, \frac{\norma{b}_{\L\infty ([0,t'']; \reali)}}{\check g}
    \, e^{2C \, t''}
    (t''-t')
    \\
    & &
    \qquad \times
    2
    \left(
      \frac{C}{\check g} + C + \hat g \, t''
      \left(\norma{d}_{\L\infty ([0,t''];\reali)}
        +
        \tv (d)
        +
        \norma{\partial_x g (\cdot, 0)}_{\L\infty ([0,t''];\reali)}
        +\tv{\partial_xg (\cdot, 0)}\right)
    \right)
    \\
    & \leq &
    \mathcal{L} \, (t''-t')
  \end{eqnarray*}
  for a suitable positive constant $\mathcal{L}$ dependent on $t'',
  \check g, \hat g, C$ and $\tv (b; [0,t''])$, $\norma{b}_{\L\infty
    ([0,t'']; \reali)}$.  Passing to~\eqref{eq:LipT2},
  use~\eqref{eq:C} and~\eqref{eq:Linfty}:
  \begin{eqnarray*}
    \int_{\gamma (t')}^{\gamma (t'')} \modulo{u (t'',x) - u (t',x)} \d{x}
    & \leq &
    2 \, \hat g \,
    \max\left\{
      \norma{u (t')}_{\L\infty (\reali^+; \reali)}
      \,,\;
      \norma{u (t'')}_{\L\infty (\reali^+; \reali)}
    \right\}
    (t''-t')
    \\
    & \leq &
    2 \, \hat g \,
    \left(
      \norma{u_o}_{\L\infty (\reali^+; \reali)}
      +
      \frac{1}{\check g} \, \norma{b}_{\L\infty ([0,t''];\reali)}
    \right)
    e^{C t''}
    (t''-t')
  \end{eqnarray*}
  Finally, deal with~\eqref{eq:LipT3} using~\eqref{eq:3}:
  \begin{eqnarray*}
    & &
    \int_{\gamma (t'')}^{+\infty} \modulo{u (t',x) - u (t'',x)} \d{x}
    \\
    & \leq &
    \int_{\gamma (t'')}^{+\infty}
    \modulo{u (t',x) - u \left(t', X (t';t'',x)\right)}
    \exp\left[
      \int_{t'}^{t''}
      \!\!\!
      \left(
        d \left(\tau, X (\tau; t', x)\right)
        -
        \partial_x g \left(\tau, X (\tau; t', x)\right)
      \right)
      \d{\tau}
    \right]
    \d{x}
    \\
    & \leq &
    \int_{\gamma (t'')}^{+\infty}
    \modulo{u (t',x) - u \left(t', X (t';t'',x)\right)}
    \d{x}
    \\
    & &
    +
    \int_{\gamma (t'')}^{+\infty}
    \modulo{u \left(t', X (t';t'',x)\right)}
    \modulo{
      \exp\left(
        \int_{t'}^{t''}
        \left(
          d \left(\tau, X (\tau; t', x)\right)
          -
          \partial_x g \left(\tau, X (\tau; t', x)\right)
        \right)
        \d{\tau}
      \right)
      -1
    }
    \d{x}
    \\
    & \leq &
    \tv (u) \, \norma{g (t)}_{\L\infty (\reali^+; \reali)} (t''-t')
    +
    \norma{u (t')}_{\L1 (\reali^+; \reali)}
    \, \left(\exp\left(2 C (t''-t')\right) -1\right)
    \\
    & \leq &
    \left(
      \hat g \, \tv (u)
      +
      2C \, \norma{u (t')}_{\L1 (\reali^+; \reali)}\right)
    (t''-t')
  \end{eqnarray*}
  Completing the proof of~\eqref{eq:lip-dependence}.
\end{proofof}

The following elementary lemma is of use below.

\begin{lemma}
  \label{lem:trivial}
  Let $H,K \in \reali^+$ and assume that the numbers $B_k \in
  \reali^+$ satisfy $B_{k+1} \leq H + K \, B_k$ for all $k \in
  \naturali$. Then, $B_k \leq K^k B_0 + \frac{1-K^k}{1-K} \, H$.
\end{lemma}

\begin{proofof}{Theorem~\ref{thm:wp}}
  Fix a time $T$ so that
  \begin{equation}
    \label{eq:2}
    \gamma (T) \in \bigl]0, \min_{i=1,\ldots,n} \bar x_i \bigr[
  \end{equation}
  and define $u^0 (t,x) = u_o (x)$ for $t\in [0,T]$. Recursively, for
  $k\geq 1$ let $u^k \equiv (u^k_1, \ldots, u^k_n)$ solve
  \begin{equation}
    \label{eq:k}
    \left\{
      \begin{array}{@{}l@{\quad}r@{\;}c@{\;}l@{}}
        \partial_t u^k_i + \partial_x \left(g_i (t,x) \, u^k_i\right)
        =
        d_i (t,x) \, u^k_i (t)
        & (t,x) & \in & \reali^+ \times \reali^+
        \\
        g_i (t, 0) \, u^k_i (t, 0)
        =
        b^k_i (t)
        & t & \in & \reali^+
        \\
        u_i (0, x) =  u_i^o (x)
        & x & \in & \reali^+
      \end{array}
    \right.
    \qquad i=1, \ldots, n
  \end{equation}
  where
  \begin{equation}
    \label{eq:bk}
    b^k_i (t)
    =
    \alpha_i\left(t,u^{k-1}_1 (t, \bar x_1), \ldots, u^{k-1}_n (t, \bar x_n)\right)
    +
    \beta_i\left(
      \int_{I_1} u^{k-1}_1 (t,x) \d{x}, \ldots,
      \int_{I_n} u^{k-1}_n (t,x) \d{x}
    \right) \,.
  \end{equation}
  Note that~\textbf{(b)} is satisfied and Lemma~\ref{lem:stability}
  applies. Indeed, if $k=1$, then $b^1_i$ is independent on time. Let
  $k >1$, then by~\eqref{eq:5} and~(\ref{eq:tvi})
  \begin{eqnarray*}
    & &
    \tv\left(b^k_i; [0,T]\right)
    \\
    & \leq &
    \tv\left(
      \alpha_i\left(\cdot,u_j^{k-1} (\cdot, \bar x)_{|j=1, \ldots,n}\right); [0,T]
    \right)
    +
    \tv \left(\beta_i\left(
        \int_{I_j} u^{k-1}_j (\cdot,x) \d{x}_{|j=1, \ldots,n}
      \right); [0,T]
    \right)
    \\
    & \leq &
    \Lip (\alpha)
    \left(
      T
      +
      \tv\left(u_j^{k-1} (\cdot, \bar x)_{|j=1, \ldots,n}; [0,T]\right)
    \right)
    +
    \Lip (\beta)
    \sum_{i=1}^n
    \tv \left(\int_{I_i} u^{k-1}_i (\cdot,x) \d{x};[0,T]\right)
    \\
    & \leq &
    \Lip (\alpha)
    \left(
      T
      +
      \tv\left(u_j^{k-1} (\cdot, \bar x_j)_{|j=1, \ldots,n}; [0,T]\right)
    \right)
    \\
    & &
    +
    C \Lip (\beta) \sum_{i=1}^n
    \int_0^t \left[
      \norma{u^{k-1}_i (\tau)}_{\L\infty (\reali^+; \reali)}
      +
      \tv\left(u^{k-1}_i (\tau)\right)
    \right]
    \d\tau
    \\
    & \leq &
    \Lip (\alpha)
    \left(
      T
      +
      \tv\left(u_j^{k-1} (\cdot, \bar x_j)_{|j=1, \ldots,n}; [0,T]\right)
    \right)
    \\
    & &
    +
    C \, \Lip (\beta)
    \int_0^t \left(
      \norma{u^{k-1} (\tau)}_{\L\infty (\reali^+; \reali^n)}
      +
      \tv\left(u^{k-1} (\tau)\right)
    \right)
    \d\tau
  \end{eqnarray*}
  which is finite by induction. Lemma~\ref{lem:stability} then ensures
  existence and uniqueness of a solution
  to~\eqref{eq:k}--\eqref{eq:bk} for any $k > 0$. By construction, the
  choice~\eqref{eq:2} ensures that
  \begin{equation}
    \label{eq:furbata}
    u^k (t,x) = u^1 (t,x)
    \quad \mbox{ for all } \quad x > \gamma (t)
    \quad \mbox{ and } \quad k \geq 1 \,.
  \end{equation}
  Therefore, also $\alpha_i\left(t, u^{k} (t, \bar x_i)\right) =
  \alpha_i\left(t, u^{1} (t, \bar x_i)\right)$ for all $t \in [0,T]$,
  for all $k \geq 1$ and all $i = 1, \ldots, n$. Compute now
  \begin{eqnarray}
    \nonumber
    & &
    \norma{u^{k+1}_i (t) - u^{k}_i (t)}_{\L1 (\reali^+;\reali)}
    \\
    \nonumber
    & \leq &
    \frac{1}{\check g} \norma{b^{k+1}_i - b^{k}_i}_{\L1 ([0,t];\reali)}
    \\
    \nonumber
    & \leq &
    \norma{
      \alpha_i\left(\cdot, u^{k} (\cdot, \bar x_i)\right) -
      \alpha_i\left(\cdot, u^{k-1} (\cdot, \bar x_i)\right)}_{\L1 ([0,t];\reali)}
    \\
    \nonumber
    & + &
    \textstyle
    \norma{
      \beta_i\!\left(
        \int_{I_1} \!\! u^{k}_1 (\cdot,x) \d{x}, \ldots,
        \int_{I_n} \!\! u^{k}_n (\cdot,x) \d{x}
      \right)
      -
      \beta_i\!\left(
        \int_{I_1} \!\! u^{k-1}_1 (\cdot,x) \d{x}, \ldots,
        \int_{I_n} \!\! u^{k-1}_n (\cdot,x) \d{x}
      \right)}_{\L1 ([0,t];\reali)}
    \\
    \nonumber
    & \leq &
    \Lip (\beta)
    \sum_{j=1}^n
    \norma{
      \int_{I_j}
      \left(u^k_j (\cdot, x) - u^{k-1}_j (\cdot,x) \right) \d{x}
    }_{\L1 ([0,t];\reali)}
    \\
    \nonumber
    & \leq &
    \Lip (\beta)
    \sum_{j=1}^n
    \norma{
      \int_{\reali^+}
      \left(u^k_j (\cdot, x) - u^{k-1}_j (\cdot,x) \right) \d{x}
    }_{\L1 ([0,t];\reali)}
    \\
    \nonumber
    & \leq &
    \Lip (\beta)
    \sum_{j=1}^n
    \int_0^t
    \norma{u^k_j (\tau) - u^{k-1}_j (\tau)}_{\L1 (\reali^+; \reali)} \d{\tau}
    \\
    \label{eq:tvb}
    & \leq &
    \Lip (\beta)
    \int_0^t
    \norma{u^k (\tau) - u^{k-1} (\tau)}_{\L1 (\reali^+; \reali^n)} \d{\tau}
  \end{eqnarray}
  Adding up all the components,
  \begin{eqnarray*}
    \norma{u^{k+1} - u^{k} }_{\C0([0,T];\L1 (\reali^+;\reali^n))}
    & \leq &
    n \, \Lip(\beta) \,
    \int_0^T \norma{u^k - u^{k-1}}_{\L1 (\reali^+; \reali^n)} \d\tau
    \\
    & \leq &
    n \, \Lip (\beta) \, T \,
    \norma{u^k - u^{k-1} }_{\C0([0,T];\L1 (\reali^+; \reali^n))}
  \end{eqnarray*}
  and, recursively,
  \begin{displaymath}
    \norma{u^{k+1} - u^{k}}_{\C0([0,T];\L1 (\reali^+;\reali^n))}
    \leq
    \left(n \, \Lip (\beta) \, T\right)^k
    \norma{u^1 (\tau) - u^0 (\tau)}_{\C0([0,T];\L1 (\reali^+; \reali^n))} .
  \end{displaymath}
  Choosing now also $T < 1 / \left(n\,\Lip (\beta)\right)$, the
  sequence $u^k$ is a Cauchy sequence and we obtain the existence of a
  map $u^* \in \C0 \left([0,T];\L1 (\reali^+; \reali^n)\right)$ which
  is the limit of the sequence $u^k$, in the sense that
  \begin{equation}
    \label{eq:limit}
    \lim_{k \to +\infty}
    \sup_{t \in [0,T]}
    \norma{u^k (t) - u^* (t)}_{\L1 (\reali^+; \reali^n)} = 0 \,.
  \end{equation}

  To prove that $u^*$ solves~\eqref{eq:1}, it is sufficient to check
  that the boundary condition is attained. Indeed, proving that $u^*$
  is a weak solution to the balance law is a standard
  procedure. Clearly, the initial datum is attained, since $u^k (0) =
  u_o$ for all $k$.

  Using Lemma~\ref{lem:stability}, \eqref{eq:12}, \eqref{eq:furbata}
  and~\eqref{eq:limit}, for all large $k \in \naturali$ we have
  \begin{eqnarray}
    \nonumber
    \!\!\!\!\!\!\!\!\!
    \norma{b^{k+1}}_{\L\infty ([0,T]; \reali^n)}
    & \leq &
    \Lip (\alpha) \, \sum_{i=1}^n
    \norma{u^{k} (\cdot, \bar x_i)}_{\L\infty ([0,T];\reali)}
    +
    \Lip (\beta) \, \norma{u^{k}}_{\C0 ([0,T]; \L1 (\reali^+; \reali^n))}
    \\
    \label{eq:normaB}
    & \leq &
    \Lip (\alpha) \,
    \norma{u_o}_{\L\infty (\reali^+;\reali^n)} e^{C T}
    +
    \Lip (\beta) \, (\norma{u^*}_{\C0 ([0,T]; \L1 (\reali^+; \reali^n))}+1) .
  \end{eqnarray}
  With the above choice of $T$, using~\eqref{eq:tvb} and
  Lemma~\ref{lem:stability},
  \begin{eqnarray*}
    \!\!\! & &
    \tv\left(b^{k+1}_i;[0,T]\right)
    \\
    \!\!\! & \leq &
    \Lip (\alpha) \, T
    +
    \Lip (\alpha)
    \left(
      \norma{u_o}_{\L\infty ([0,\max_j \bar x_j]; \reali^n)}
      +
      \tv\left(u_o; [0,\max_j \bar x_j]\right)
    \right)  e^{C T}
    \\
    \!\!\! & &
    +
    C \, \Lip (\beta)
    \int_0^T
    \left(
      \norma{u_o}_{\L\infty (\reali^+; \reali^n)}
      +
      \frac{1}{\check g} \norma{b^k}_{\L\infty ([0,\tau]; \reali^n)}
    \right)
    e^{C \tau}
    \d\tau
    \\
    \!\!\! & &
    +
    C \Lip (\beta) \!\!
    \int_0^T
    \!\!
    \left[
      \norma{u_o}_{\L\infty (\reali^+; \reali^n)}
      +
      \tv (u_o)
      +
      \frac{1}{\check g}
      \left[
        \frac{C+\check g}{\check g}
        \norma{b^k}_{\L\infty ([0,\tau]; \reali^n)}
        \!\!\! \! {+}
        \tv (b^k;[0,\tau])
      \right]
    \right]\!
    e^{C \tau}
    \d\tau
    \\
    \!\!\! & \leq &
    \Lip (\alpha) \, T
    +
    \Lip (\alpha) \, \left(\norma{u_o}_{\L\infty (\reali^+; \reali^n)}
      + \tv (u_o) \right)
    e^{CT}
    \\
    \!\!\! & &
    +
    2 \, C \, \Lip (\beta)
    \int_0^T
    \left(
      \norma{u_o}_{\L\infty (\reali^+; \reali^n)}
      +
      \tv (u_o)
    \right)
    e^{C \tau}
    \d\tau
    \\
    \!\!\! & &
    +
    \frac{C}{\check g} \, \Lip (\beta)
    \int_0^T
    \left(
      \frac{2\check g + C}{\check g} \,
      \norma{b^k}_{\L\infty ([0,\tau]; \reali^n)}
      +
      \tv (b^k;[0,\tau])
    \right)
    e^{C \tau}
    \d\tau
    \\
    \!\!\! & \leq &
    \Lip (\alpha) \, T
    +
    \Lip (\alpha) \, \left(\norma{u_o}_{\L\infty (\reali^+; \reali^n)}
      + \tv (u_o) \right) e^{CT}
    \\
    \!\!\! & &
    +
    2 \, \Lip (\beta)
    \left(
      \norma{u_o}_{\L\infty (\reali^+; \reali^n)}
      +
      \tv (u_o)
    \right)
    (e^{C T} -1)
    \\
    \!\!\! & &
    +
    \frac{C}{\check g} \, \Lip (\beta)
    \int_0^T
    \left(
      \frac{2\check g + C}{\check g} \,
      \norma{b^k}_{\L\infty ([0,\tau]; \reali^n)}
      +
      \tv (b^k;[0,\tau])
    \right)
    e^{C \tau}
    \d\tau
    \\
    \!\!\! & \leq &
    \Lip (\alpha) \, T
    +
    \Lip (\alpha) \, \left(\norma{u_o}_{\L\infty (\reali^+; \reali^n)}
      + \tv (u_o) \right) e^{CT}
    \\
    \!\!\! & &
    +
    2 \, \Lip (\beta)
    \left(
      \norma{u_o}_{\L\infty (\reali^+; \reali^n)}
      +
      \tv (u_o)
    \right)
    (e^{C T} -1)
    \\
    \!\!\! & &
    +
    \frac{2\check g + C}{\check g^2} \, \Lip (\beta)
    \left(
      \norma{b^k}_{\L\infty ([0,T]; \reali^n)}
      +
      \tv (b^k;[0,T])
    \right)
    (e^{C T}-1)
  \end{eqnarray*}
  Inserting now the estimate~\eqref{eq:normaB} in the latter term
  above, we can apply Lemma~\ref{lem:trivial} to the inequality
  $B_{k+1} \leq H + K \, B_k$, where
  \begin{eqnarray*}
    B_k  & = & \tv\left(b^k;[0,T]\right)
    \\
    H & = &
    n \, \Lip (\alpha) \, T
    +
    n
    \left(
      \Lip (\alpha) \, e^{CT}
      +
      2 \, \Lip (\beta)   (e^{C T} -1)
    \right)
    \left(
      \norma{u_o}_{\L\infty (\reali^+; \reali^n)}
      + \tv (u_o)
    \right)
    \\
    & &
    +
    n \, \frac{2\check g + C}{\check g^2} \, \Lip (\beta) \, (e^{C T} -1)
    \\
    & &
    \qquad \times
    \left(
      \Lip (\alpha) \,
      \norma{u_o}_{\L\infty (\reali^+;\reali^n)} e^{C T}
      +
      \Lip (\beta) \, (\norma{u^*}_{\C0 ([0,T]; \L1 (\reali^+; \reali^n))}+1)
    \right)
    (e^{C T}-1)
    \\
    K & = &
    n \, \frac{2\check g + C}{\check g^2} \, \Lip (\beta) \, (e^{C T}-1) \,.
  \end{eqnarray*}
  Hence, as soon as $T$ is so small that
  \begin{displaymath}
    n \, \frac{2\check g + C}{\check g^2} \, \Lip (\beta) \, (e^{C T}-1) < 1 \,,
  \end{displaymath}
  we obtain a bound on $\tv\left(b^k;[0,T]\right)$ uniform in
  $k$. This bound, thanks to Lemma~\ref{lem:stability}, ensures that
  also $\sup_{t \in [0,T]} \sup _{k \in \naturali} \tv\left(u^k
    (t)\right) < +\infty$ and, by the lower semicontinuity of the
  total variation with respect to the $\L1$ topology, also $\sup_{t
    \in [0,T]} \tv\left(u^* (t)\right) < +\infty$. Therefore, the
  trace $\lim_{x \to 0+} u^*(t,x)$ exists for all $t \in [0,T]$.

  The uniform bound on $\tv\left(b^k;[0,T]\right)$, together
  with~\eqref{eq:lip-dependence}
  and~\cite[Theorem~2.4]{BressanLectureNotes}, ensure that for a.e.~$x
  \in \reali^+$, we have $\lim_{k\to+\infty} u^k (t,x) = u^*
  (t,x)$. Choose one such $x$ and observe that:
  \begin{eqnarray*}
    \!\! u^* (t,x) \!\!
    & = &
    \lim_{k\to+\infty} u^k (t,x)
    \\
    & = &
    \lim_{k\to+\infty}
    \frac{b^k\left(T (0;t,x)\right)}{g\left(T (0;t,x),0\right)}
    \exp \left(
      \int_{T (0;t,x)}^t
      \left(
        d(\tau,X (\tau;t,x))
        -
        \partial_x g \left(\tau,X (\tau;t,x)\right)
      \right)
      \d\tau
    \right)
    \\
    & = &
    \frac{\lim_{k\to+\infty}b^k\left(T (0;t,x)\right)}{g\left(T (0;t,x),0\right)}
    \exp \left(
      \int_{T (0;t,x)}^t
      \left(
        d(\tau,X (\tau;t,x))
        -
        \partial_x g \left(\tau,X (\tau;t,x)\right)
      \right)
      \d\tau
    \right)
    \\
    & = &
    \frac{\lim_{k\to+\infty}
      \mathcal{B} \!\left(T (0;t,x), u^{k-1}\right)}{g\left(T (0;t,x),0\right)}
    \exp \left[
      \int_{T (0;t,x)}^t
      \!\!\!
      \left(
        d(\tau,X (\tau;t,x))
        -
        \partial_x g \left(\tau,X (\tau;t,x)\right)
      \right)
      \d\tau
    \right]
    \\
    & = &
    \frac{\mathcal{B} \!\left(T (0;t,x), u\right)}{g\left(T (0;t,x),0\right)}
    \exp \left(
      \int_{T (0;t,x)}^t
      \left(
        d(\tau,X (\tau;t,x))
        -
        \partial_x g \left(\tau,X (\tau;t,x)\right)
      \right)
      \d\tau
    \right)
  \end{eqnarray*}
  where in the last step above we used the convergences:
  \begin{displaymath}
    \begin{array}{rcll}
      \displaystyle
      \lim_{k\to+\infty} \alpha_i\left(t, u_j^k (t,\bar x_j)_{|j=1,\ldots,n}\right)
      & = &
      \displaystyle
      \alpha_i\left(t, u^* (t, \bar x_j)_{|j=1,\ldots,n}\right)
      & \mbox{by~\eqref{eq:furbata}.}
      \\
      \displaystyle
      \lim_{k\to+\infty} \beta_i\left(\int_{I_j} u_j^k (t,x) \d{x}_{|j=1,\ldots,n}\right)
      & = &
      \displaystyle
      \beta_i\left(\int_{I_j} u^* (t,x) \d{x}_{|j=1,\ldots,n}\right)
      & \mbox{by~\eqref{eq:limit}.}
    \end{array}
  \end{displaymath}

  The time $T$ chosen above depends only on $\beta$, $\min_i \bar
  x_i$, on $d$ and on $g$. In particular, it is independent from the
  initial datum. Hence, the above procedure can be iterated, extending
  $u^*$ to a function defined on all $\reali^+$, i.e.~$u^* \in
  \C0\left(\reali^+; \L1 (\reali^+; \reali^n)\right)$.

  \smallskip

  Let now $u_o', u_o''$ be two initial data. Define $\bar x =
  \min_{i=1, \ldots, n} \bar x_i$ and $\bar t = \Gamma (\bar
  x)$. Denote $\bar I = \left[\bar x, +\infty\right[$. To prove the
  stability estimate, with obvious notation, preliminarily compute for
  $t \in [0,\bar t]$:
  \begin{eqnarray*}
    \norma{u'' (t) - u' (t)}_{\L1 (\bar I; \reali^n)}
    & \leq &
    \norma{u_o''- u_o'}_{\L1 (\reali^+; \reali^n)} e^{C t}
    \\
    \norma{u'' (t) - u' (t)}_{\L\infty (\bar I; \reali^n)}
    & \leq &
    \norma{u_o''- u_o'}_{\L\infty (\reali^+; \reali^n)} e^{C t}
  \end{eqnarray*}
  Moreover, by~\eqref{eq:bk} and~\eqref{eq:5}, for $t \in [0, \bar
  t]$,
  \begin{eqnarray}
    \nonumber
    \norma{b'' - b'}_{\L\infty ([0,\bar t];\reali^n)}
    & \leq &
    n\, \Lip (\alpha) \, \norma{u'' (t) - u' (t)}_{\L\infty (\bar I; \reali^n)}
    +
    n \, \Lip(\beta) \, \norma{u'' (t) - u' (t)}_{\L1 (\reali^+; \reali^n)}
    \\
    \label{eq:bInfty}
    & \leq &
    n\, \Lip (\alpha) \,
    \norma{u_o''- u_o'}_{\L\infty (\reali^+; \reali^n)} e^{C t}
    \\
    \nonumber
    & &
    +
    n \, \Lip(\beta) \, \norma{u'' (t) - u' (t)}_{\L1 (\reali^+; \reali^n)} \,.
    \\
    \nonumber
    \norma{b'' - b'}_{\L1 ([0,\bar t];\reali^n)}
    & \leq &
    n\, \Lip (\alpha) \, \frac{e^{C t}-1}{C}
    \norma{u_o''- u_o'}_{\L\infty (\reali^+; \reali^n)}
    \\
    \label{eq:b1}
    & &
    +
    n \, \Lip(\beta)
    \int_0^t \norma{u'' (\tau) - u' (\tau)}_{\L1 (\reali^+; \reali^n)} \d\tau\,.
  \end{eqnarray}
  Therefore, by~\eqref{eq:lip-dependence}, for $t \in [0, \bar t]$,
  \begin{eqnarray*}
    \norma{u'' (t) - u' (t)}_{\L1 (\reali; \reali^n)}
    & \leq &
    \left(
      \frac{1}{\check g} \, \norma{b''-b'}_{\L1 ([0,t]; \reali)}
      +
      \norma{u_o''- u_o'}_{\L1 (\reali^+; \reali^n)}
    \right)
    e^{C t}
    \\
    & \leq &
    \frac{n\, \Lip (\alpha)}{\check g} \,
    \, \frac{e^{C t}-1}{C} \, e^{C t}
    \norma{u_o''- u_o'}_{\L\infty (\reali^+; \reali^n)}
    \\
    & &
    +
    \norma{u_o''- u_o'}_{\L1 (\reali^+; \reali^n)} \, e^{C t}
    \\
    & &
    +
    n \, \Lip(\beta)
    \int_0^t \norma{u'' (\tau) - u' (\tau)}_{\L1 (\reali^+; \reali^n)} \d\tau
    e^{C t} \,.
  \end{eqnarray*}
  An application of Gronwall Lemma yields, for $t \in [0, \bar t]$,
  \begin{eqnarray*}
    \norma{u'' (t) - u' (t)}_{\L1 (\reali; \reali^n)}
    & \leq &
    \left(
      \frac{n \Lip (\alpha)}{\check g}
      \frac{e^{C t}-1}{C}
      \norma{u_o''- u_o'}_{\L\infty (\reali^+; \reali^n)}
      +
      \norma{u_o''- u_o'}_{\L1 (\reali^+; \reali^n)}
    \right)
    \\
    & &
    \times
    \exp\left(
      C t +
      n \Lip(\beta)
      e^{C t}
    \right)
  \end{eqnarray*}
  To iterate beyond time $\bar t$, using~\eqref{eq:Linfty},
  \eqref{eq:bInfty} and the above bound to estimate
  \begin{eqnarray*}
    & &
    \norma{u'' (t) - u' (t)}_{\L\infty (\reali; \reali^n)}
    \\
    & \leq &
    \left(
      \frac{1}{\check g} \, \norma{b''-b'}_{\L\infty ([0,t];\reali^n)}
      +
      \norma{u_o''-u_o'}_{\L\infty (\reali^+; \reali^n)}
    \right)
    e^{C t}
    \\
    & \leq &
    \frac{n\, \Lip (\alpha)}{\check g}
    \norma{u_o''- u_o'}_{\L\infty (\reali^+; \reali^n)} e^{2C t}
    +
    \frac{n\, \Lip (\beta)}{\check g}
    \norma{u'' (t) - u' (t)}_{\L1 (\reali^+; \reali^n)} e^{C t}
    \\
    & &
    +
    \norma{u_o''-u_o'}_{\L\infty (\reali^+; \reali^n)}
    e^{C t}
    \\
    & \leq &
    \left(
      \left(\frac{n\, \Lip (\alpha)}{\check g} \, e^{C t} + 1\right)
      \norma{u_o''- u_o'}_{\L\infty (\reali^+; \reali^n)}
      +
      \frac{n\, \Lip (\beta)}{\check g}
      \norma{u'' (t) - u' (t)}_{\L1 (\reali^+; \reali^n)}
    \right)
    e^{C t}
    \\
    & \leq &
    \left(
      \left(\frac{n\, \Lip (\alpha)}{\check g} \, e^{C t} + 1\right)
      +
      \frac{n\, \Lip (\beta)}{\check g}
      \frac{n \Lip (\alpha)}{\check g}
      \frac{e^{C t}-1}{C}
      \exp\left(
        C t
        +
        n \Lip(\beta)
        e^{C t}
      \right)
    \right)
    \\
    & &
    \qquad
    \times
    \norma{u_o''- u_o'}_{\L\infty (\reali^+; \reali^n)}
    e^{C t}
    \\
    & &
    +
    \frac{n\, \Lip (\beta)}{\check g}
    \exp\left(
      2C t
      +
      n \Lip(\beta)
      e^{C t}
    \right)
    \norma{u_o'' - u_o'}_{\L1 (\reali^+; \reali^n)}
  \end{eqnarray*}
  completing the proof.
\end{proofof}

\begin{proofof}{Theorem~\ref{thm:Stab}}
  Define, for $i=1, \ldots, n$,
  \begin{eqnarray*}
    b_i' (t)
    & = &
    \alpha_i'\left(t, u_1' (t,\bar x_1-), \ldots, u_n' (t, \bar x_n-)\right)
    +
    \beta_i'\left(
      \int_{I_1} u_1' (t,x) \d{x}, \ldots, \int_{I_n} u_n' (t,x) \d{x}
    \right),
    \\
    b_i'' (t)
    & = &
    \alpha_i''\left(t, u_1'' (t,\bar x_1-), \ldots, u_n'' (t, \bar x_n-)\right)
    +
    \beta_i''\left(
      \int_{I_1} u_1'' (t,x) \d{x}, \ldots, \int_{I_n} u_n'' (t,x) \d{x}
    \right).
  \end{eqnarray*}
  Preliminary, using~\eqref{eq:5}, let us estimate the term
  \begin{eqnarray*}
    \norma{b'-b''}_{\L1 ([0,t]; \reali^n)} \!\!\!
    & \leq &
    \!\!\!
    \sum_{i=1}^n
    \int_0^t \modulo{\alpha'_i\left(s, u'_1(s, \bar x_1), \ldots,
        u'_n(s, \bar x_n)\right) - \alpha''_i\left(s,u''_1(s, \bar
        x_1), \ldots, u''_n(s, \bar x_n)\right)} \d s
    \\
    & &
    + \sum_{i=1}^n \int_0^t \modulo{\beta'_i\left(\int_{I_j} u'_j(s,
        x) \d x \right) - \beta''_i\left(\int_{I_j} u''_j(s, x) \d
        x\right)} \d s
    \\
    & \leq &
    \!\!\!
    t \norma{\alpha' -
      \alpha''}_{\C0\left(\reali^+\times\reali^n; \reali^n\right)} +
    \Lip\left(\alpha''\right) \sum_{i=1}^n \int_0^t \modulo{u_j'(s,
      \bar x_j) - u_j''(s, \bar x_j)} \d s
    \\
    & &
    + t \norma{\beta' - \beta''}_{\C0\left(\reali^n;
        \reali^n\right)} + \Lip\left(\beta''\right) \int_0^t
    \norma{u'(s) - u''(s)}_{\L1\left(\reali^+; \reali^n\right)} \d s \,.
  \end{eqnarray*}
  Define $\bar x = \min_{i=1, \ldots, n} \bar x_i$ and $\bar t =
  \Gamma (\bar x)$. As long as $s \in [0, \bar T]$, we have $u_j'(s,
  \bar x_j) = u_j''(s, \bar x_j)$. Hence the above estimate leads to
  \begin{eqnarray*}
    \norma{b'-b''}_{\L1 ([0,t]; \reali^n)}
    & \leq &
    t
    \left(
      \norma{\alpha' - \alpha''}_{\C0\left(\reali^+\times\reali^n; \reali^n\right)}
      +
      \norma{\beta' - \beta''}_{\C0\left(\reali^n; \reali^n\right)}
    \right)
    \\
    & &
    +
    \Lip\left(\beta''\right) \int_0^t
    \norma{u'(s) - u''(s)}_{\L1\left(\reali^+; \reali^n\right)} \d s
  \end{eqnarray*}
  for all $t \in [0,\bar t]$. In the same time interval,
  \begin{eqnarray*}
    \norma{u' (t) - u'' (t)}_{\L1 (\reali^+; \reali^n)}
    & \leq &
    \frac{e^{C t}}{\check g} \, \norma{b'-b''}_{\L1 ([0,t]; \reali^n)}
    \\
    & \leq &
    \frac{t \, e^{C t}}{\check g}
    \left(
      \norma{\alpha' - \alpha''}_{\C0\left(\reali^+\times\reali^n; \reali^n\right)}
      +
      \norma{\beta' - \beta''}_{\C0\left(\reali^n; \reali^n\right)}
    \right)
    \\
    & &
    +
    \frac{\Lip(\beta'') \, e^{C t}}{\check g} \int_0^t
    \norma{u'(s) - u''(s)}_{\L1\left(\reali^+; \reali^n\right)} \d s
  \end{eqnarray*}
  so that by Gronwall Lemma, for $t \in [0, \bar t]$,
  \begin{eqnarray*}
    \norma{u' (t) - u'' (t)}_{\L1 (\reali^+; \reali^n)}
    & \leq &
    \left[
      \frac{t \, e^{C t}}{\check g}
      \left(
        \norma{\alpha' - \alpha''}_{\C0\left(\reali^+\times\reali^n; \reali^n\right)}
        +
        \norma{\beta' - \beta''}_{\C0\left(\reali^n; \reali^n\right)}
      \right)
    \right]
    \\
    & &
    \times
    \exp\left[
      \frac{\Lip(\beta'') \, t \, e^{C t}}{\check g}
    \right] \,.
  \end{eqnarray*}
  A repeated application of the estimate above on the intervals
  $[(k-1)\bar t, k\, \bar t]$ allows to complete the proof.
\end{proofof}

\begin{proofof}{Corollary~\ref{cor:Azmy}}
  Introduce $u_1,u_2, \ldots$ as in table~\eqref{eq:AzmyTable}. Then,
  extend $d_1, d_2$ and $g_2$ to $\reali^+ \times\reali$ maintaining
  the required regularity and bounds on the total variation. The
  resulting system fits into~\eqref{eq:1}--\eqref{eq:B}. Hence,
  \textbf{(d)}, \textbf{(g)} and~\eqref{eq:5}
  hold. Theorem~\ref{thm:wp} applies, ensuring the well posedness of
  the Cauchy problem. Finally, the solution to~\eqref{eq:Azmy} is
  obtained restricting the solution
  to~\eqref{eq:1}--\eqref{eq:B}--\eqref{eq:AzmyTable} to $[0,
  a_{\max}]$ and to $[x_{\min}, x_{\max}]$.
\end{proofof}

\smallskip

\noindent\textbf{Acknowledgment:} The authors were supported by the INDAM--GNAMPA project
\emph{Conservation Laws: Theory and Applications}.

{\small

  \bibliographystyle{abbrv}

  \bibliography{Colombo-Garavello}

\begin{thebibliography}{10}

\bibitem{Ackleh2009}
A.~S. Ackleh and K.~Deng.
\newblock A nonautonomous juvenile-adult model: well-posedness and long-time
  behavior via a comparison principle.
\newblock {\em SIAM J. Appl. Math.}, 69(6):1644--1661, 2009.

\bibitem{Ackleh2012}
A.~S. Ackleh, K.~Deng, and X.~Yang.
\newblock Sensitivity analysis for a structured juvenile--adult model.
\newblock {\em Comput. Math. Appl.}, 64(3):190--200, 2012.

\bibitem{AmbrosioFuscoPallara}
L.~Ambrosio, N.~Fusco, and D.~Pallara.
\newblock {\em Functions of bounded variation and free discontinuity problems}.
\newblock Oxford Mathematical Monographs. The Clarendon Press Oxford University
  Press, New York, 2000.

\bibitem{BardosLerouxNedelec}
C.~Bardos, A.~Y. le~Roux, and J.-C. N{\'e}d{\'e}lec.
\newblock First order quasilinear equations with boundary conditions.
\newblock {\em Comm. Partial Differential Equations}, 4(9):1017--1034, 1979.

\bibitem{BrauerCastillo-Chavez}
F.~Brauer and C.~Castillo-Chavez.
\newblock {\em Mathematical models in population biology and epidemiology},
  volume~40 of {\em Texts in Applied Mathematics}.
\newblock Springer, New York, second edition, 2012.

\bibitem{BressanLectureNotes}
A.~Bressan.
\newblock {\em Hyperbolic systems of conservation laws}, volume~20 of {\em
  Oxford Lecture Series in Mathematics and its Applications}.
\newblock Oxford University Press, Oxford, 2000.
\newblock The one-dimensional Cauchy problem.

\bibitem{ColomboGuerra2010}
R.~M. Colombo and G.~Guerra.
\newblock On general balance laws with boundary.
\newblock {\em J. Differential Equations}, 248(5):1017--1043, 2010.

\bibitem{ColomboGuerraHertySchleper}
R.~M. Colombo, G.~Guerra, M.~Herty, and V.~Schleper.
\newblock Optimal control in networks of pipes and canals.
\newblock {\em SIAM J. Control Optim.}, 48(3):2032--2050, 2009.

\bibitem{GaravelloPiccoliBook}
M.~Garavello and B.~Piccoli.
\newblock {\em Traffic flow on networks}, volume~1 of {\em AIMS Series on
  Applied Mathematics}.
\newblock American Institute of Mathematical Sciences (AIMS), Springfield, MO,
  2006.
\newblock Conservation laws models.

\bibitem{Keyfitz1972}
N.~Keyfitz.
\newblock The mathematics of sex and marriage.
\newblock In {\em Proceedings of the Sixth Berkeley Symposium on Mathematical
  Statistics and Probability, Volume 4: Biology and Health}, pages 89--108,
  Berkeley, Calif., 1972. University of California Press.

\bibitem{Kruzkov}
S.~N. Kru{\v{z}}hkov.
\newblock First order quasilinear equations with several independent variables.
\newblock {\em Mat. Sb. (N.S.)}, 81 (123):228--255, 1970.

\bibitem{LeVequeBook2002}
R.~J. LeVeque.
\newblock {\em Finite volume methods for hyperbolic problems}.
\newblock Cambridge Texts in Applied Mathematics. Cambridge University Press,
  Cambridge, 2002.

\bibitem{Manfredi1993}
P.~Manfredi.
\newblock Logistic effects in the two-sex model with "harmonic mean" fertility
  function.
\newblock {\em Genus}, 49(1-2):43--65, 1993.

\bibitem{PerthameBook}
B.~Perthame.
\newblock {\em Transport equations in biology}.
\newblock Frontiers in Mathematics. Birkh\"auser Verlag, Basel, 2007.

\bibitem{Schoen1981}
R.~Schoen.
\newblock The harmonic mean as the basis of a realistic two-sex marriage model.
\newblock {\em Demography}, 18(2):201--216, 1981.
\newblock cited By (since 1996)20.

\bibitem{Schoen1983}
R.~Schoen.
\newblock Relationships in a simple harmonic mean two-sex fertility model.
\newblock {\em Journal of Mathematical Biology}, 18(3):201--211, 1983.

\bibitem{Schoen1988}
R.~Schoen.
\newblock {\em Modeling multigroup populations}.
\newblock Springer, 1988.

\bibitem{SerreII}
D.~Serre.
\newblock {\em Systems of conservation laws. 2}.
\newblock Cambridge University Press, Cambridge, 2000.
\newblock Geometric structures, oscillations, and initial-boundary value
  problems, Translated from the 1996 French original by I. N. Sneddon.

\bibitem{SundelofAberg2006}
A.~Sundelof and P.~Aberg.
\newblock Birth functions in stage structured two-sex models.
\newblock {\em Ecological Modeling}, 193:787--795, 2006.

\end{thebibliography}

}

\end{document}